\documentclass[%
times
,preprint
,1p
,authoryear
]{elsarticle}

\usepackage[utf8]{inputenc}
\usepackage{amsfonts,amssymb,amsmath,amsthm,stmaryrd,enumerate}
\usepackage{lmodern,bbold}  
\usepackage{mathtools,etoolbox,xcolor}
\usepackage[mathcal]{euscript}

\usepackage[citecolor=teal,linkcolor=teal,
    urlcolor=teal,
    colorlinks]{hyperref}

\newtheorem{theorem}{Theorem}[section]
\newtheorem{lemma}[theorem]{Lemma}
\newtheorem{corollary}[theorem]{Corollary}
\newtheorem{proposition}[theorem]{Proposition}
\theoremstyle{definition}
\newtheorem*{remark}{Remark}

\newcommand{\mesh}{\operatorname{mesh}}

\DeclarePairedDelimiter{\abs}{\lvert}{\rvert}
\DeclarePairedDelimiter{\zjel}{(}{)}
\DeclarePairedDelimiter\br\lbrack\rbrack
\DeclarePairedDelimiter\event\lbrace\rbrace
\DeclarePairedDelimiterX\set[2]\{\}{%
  #1\nonscript\::\allowbreak\nonscript\:\mathopen{}#2%
}

\def\toby#1#2{\,{\buildrel #2\over#1}\,}

\newcommand\pto[1][\to]{\toby{#1}{p}}

\let\smallset\event
\let\cbr\event

\let\eps\varepsilon
\let\phi\varphi
\let\theta\vartheta

\makeatletter
  \def\isparam{\@ifnextchar{\bgroup}}%
\makeatother

\let\texexp\exp
\def\exp{\texexp\isparam{\cbr*}{}}

\let\PEfont\mathbb
\providecommand\given{}

\def\redefgiven#1{%
  \renewcommand\given{%
    \nonscript\:%
    #1\vert%
    \allowbreak%
    \nonscript\:%
    \mathopen{}%
  }%
}
\DeclarePairedDelimiterX\PEzjel[1](){\redefgiven{\delimsize}#1}
\newcommand{\PE}[1][]{%
  \PEfont{\Pe}%
  \ifblank{#1}{}{_{#1}}%
  \isparam{\PEzjel*}{}%
}
\newcommand{\E}{\def\Pe{E}\PE}
\renewcommand{\P}{\def\Pe{P}\PE}

\def\expandargi#1\relax#2{#2{#1}}
\def\expandarg#1#2{%
  \expandafter\expandargi#2\relax{#1}}
\def\firstargi#1#2\relax{#1}
\def\firstarg#1{%
  \firstargi#1\relax
}
\edef\bschar{\expandarg{\firstarg}{\string\ }}
\def\stripbsi#1#2\relax#3{%
  \ifdefstring\bschar{#1}{%
    \def#3{#2}%
  }{%
    \def#3{#1#2}%
  }%
}
\def\stripbs#1#2{%
  \stripbsi#1\relax#2%
}
\def\defname#1#2#3#4{%
  \expandarg\stripbs{\string#1}\nn
  \edef\nn{#3\nn#4}%
  \expandarg\csdef\nn{#2}%
}

\def\defcx{%
  \def\do##1{\defname{##1}{{\mathcal{##1}}}{c}{}}%
  \dolist
}
\def\deftx{%
  \def\do##1{\defname{##1}{{\tilde{##1}}}{t}{}}%
  \dolist
}
\DeclareListParser{\dolist}{}
\deftx{ABXQYZmM\mu\F\beta\tbeta\tau\Omega}
\defcx{FGAHSCDLPB}
\let\F\cF
\let\G\cG
\let\B\cB
\let\arctg\arctan
\let\tg\tan
\def\cadlag/{c\`adl\`ag}
\def\ito/{It\^o}
\def\sign{\operatorname{sign}}
\def\real{\mathbb{R}}
\def\Z{\mathbb{Z}}
\def\Q{\mathbb{Q}}
\def\interior{\operatorname{int}}

\def\Isymb{{\mathbb 1}}
\def\substr#1{_{\zjel{#1}}}
\newcommand{\I}[1][]{%
  \Isymb\ifblank{#1}{\substr}{_{#1}}%
}

\begin{document}

\begin{frontmatter}
  \title{On the lack of semimartingale property\tnoteref{t1}}
  \author[1]{Vilmos Prokaj}
  \ead{vilmos.prokaj@ttk.elte.hu}
  \author[2,1]{László Bondici \texorpdfstring{\fnref{fn2}}{}}
  \ead{bondici@eotvos.elte.hu}
  \address[1]{ELTE Eötvös Loránd University, Budapest, Hungary.}
  \address[2]{Alfréd Rényi Institute of Mathematics, Budapest, Hungary.}
  \fntext[fn2]{The majority of the work was carried out while the second author 
  was at the \texorpdfstring{Eötvös Loránd}{} University.}

  \tnotetext[t1]{The project was supported by the European Union,
    co-financed by the European Social Fund (EFOP-3.6.3-VEKOP-16-2017-00002).}

  \begin{abstract}
    In this work we extend the characterization of semimartingale functions 
    in \cite{MR597337} to the non-Markovian setting. We prove that if a function of a 
    semimartingale remains a semimartingale, then under certain conditions the function 
    must have intervals where it is a difference of two convex functions. 
    Under suitable conditions this property also holds for random functions.
    As an application, we prove that the median process defined in \cite{drift:2009} 
    is not a semimartingale. The same process appears also in \cite{MR1741809}
    where the question of the semimartingale property is raised but not settled.
  \end{abstract}
  \begin{keyword}
    semimartingale property\sep semimartingale function
    \MSC
    60H05 
    60J65\sep 
    60J55\sep 
  \end{keyword}
\end{frontmatter}

\section{Introduction}

Let $B$ be a Brownian motion and suppose that
$(D_t(x))_{t\geq0,x\in[0,1]}$ satisfies 
the stochastic differential equation
\begin{equation}\label{eq:D}
  dD_t(x)=D_t(x)\wedge (1-D_t(x)) d B_t=\sigma (D_t (x))dB_t, \quad D_0(x)=x.
\end{equation}
This two parameter process was analyzed in \cite{drift:2009} in
detail and played an important role in the construction that led to
the solution of the drift hiding problem. \cite{MR1741809} considers
\begin{equation}\label{eq:G}
  dG_t (x)=dB_t+\beta\sign (G _t(x))d t,\quad G_0 (x)=x,\quad x\in\real.
\end{equation}
When $\beta=1/2$ then $G$ is 
a transformed version of $D$, $G (p(x))=p (D (x))$, where $p$
is the Lamperti transformation associated to $\sigma$, that is
$p'=1/\sigma$, $p (1/2)=0$. The  case $\beta>0$ is just the matter of
scaling both time and space. The aim of this paper is to show that
$(D_t^{-1} (1/2))_{t\geq0}$ is not a semimartingale.
As $(G_t^{-1} (0))_{t\geq0}$ is the same as $p^{-1}\circ D^{-1}(1/2)$,
this result also confirms the expectation of \citeauthor[see the remark
after Proposition 1.2 on pp. 288]{MR1741809}. Their motivation came
from \cite{MR1681142}; $G_t (x)-B_t$ is the special case of the 
bifurcation model of that paper.

Our argument uses the fact that both $D$ and $G$ can be viewed as  a stochastic flow
in the sense of \cite{MR867686}.
For $t\geq 0$ denote
$B_{s,t}=B_{s+t}-B_{s}$, which is a Brownian motion 
that starts evolving at time $s$.  Let $D_{s,t} (x)$ be the solution of
\eqref{eq:D} starting from $x$ and driven by $(B_{s,t})_{t\geq 0}$. Then 
$D_{s+t}=D_{s,t}\circ D_s$ and $D_{s+t}^{-1} 
(1/2)=D_s^{-1} (D_{s,t}^{-1} (1/2))$ for any $s,t\geq0$.

Suppose that $m_t=D^{-1}_t (1/2)$, $t\geq0$ is a semimartingale in the
filtration of $B$ denoted by $(\F_t)_{t\geq0}$.
Fix a positive time point $s$, 
then the same is true for $(m_{s,t}=D_{s,t}^{-1} (1/2))_{t\geq 0}$ in the filtration of
$(B_{s,t})_{t\geq 0}$, denoted by $(\F_{s,t})_{t\geq0}$.
Put $\G_t=\F_s\vee\F_{s,t}$ for all $t$, then $\G$ is an initial
enlargement of $(\F_{s,t})_{t\geq0}$ with an independent $\sigma$-algebra
$\F_s$, so $(m_{s,t})_{t\geq0}$ remains a semimartingale in $\G$.

On the other hand if $m$ is a semimartingale in $\F$, then since
$\G_t=\F_{s+t}$, the process $m_{s+t}=D_{s}^{-1} (m_{s,t})$, $t\geq0$ is also a 
semimartingale in $\G$, which is obtained by substituting a $\G$-semimartingale 
into a random $\G_0$ measurable function. By the main result of 
Section \ref{sec:rand F} it means that $D^{-1}_s$ must be a semimartingale function for the
process $(m_{s,t})_{t\geq0}$ almost surely.
We give a necessary condition for this in Section \ref{sec:semi f}, 
and analysing $D_s$ and $m_{s,t}$ 
in detail in Section \ref{sec:m not semi} we conclude that $D_s^{-1}$ is almost surely 
not a semimartingale function for $(m_{s,t})_{t\geq0}$.

In Section \ref{sec:examples} we give a simple proof in a similar
fashion for a result of \cite{MR1118450}. 
Throughout the paper all $\sigma$-algebras are augmented with null sets 
and each filtration is assumed to be right continuous.

\section{Semimartingale functions}\label{sec:semi f}

We call a function $F$ a semimartingale function for $X$ if $F (X)$ is
a semimartingale. Semimartingale functions are characterized in the
Markovian setting by \cite{MR597337}. First, we give a simple proof of
their result for Brownian semimartingale functions using the embedded
random walk of the Brownian motion. This argument, with some slight
modification, provides a necessary condition for $F$ being a semimartingale
function for a continuous semimartingale $X$ even for the case when
$F$ is random.

A central part of the argument is a characterization 
of differences of convex functions. We prove below in Proposition
\ref{prop:tot var} that if  
$\limsup_{\delta\to0}\mu_{\delta,F} (I)<\infty$ for an interval $I$, then $F$ is a
difference of two convex functions on $I$, where
\begin{displaymath}
  \mu_{\delta,F}(H)=
  \sum_{k:k\delta \in H}\frac1\delta\abs*{(\Delta^\delta F)(k\delta)},
\end{displaymath}
and $(\Delta^\delta F)(x) = F(x+\delta)+F(x-\delta)-2F(x)$ is a kind
of discrete Laplace operator. To motivate this condition, observe that 
when $F$ is twice continuously differentiable then the limit is the total variation
of $F'$ on $H$.
Since in the argument we consider a single $F$, we also use the
notation $\mu_\delta$ for the measure $\mu_{\delta,F}$.

\begin{theorem}\label{thm:BM}
  Suppose $B$ is a Brownian motion and $F:\real\to\real$ is a
  continuous function.
  Then $F (B)$ is a semimartingale if and
  only if $F$ is a difference of two convex functions.
\end{theorem}

\begin{proof}
  The sufficiency of the condition in the statement is just the
  \ito/-Tanaka formula. So we only show the necessity.

  In the proof we use the embedded random walk of the Brownian motion
  $B$ with step size $\delta$. That is, we define the stopping times
  \begin{displaymath}
    \tau_0^\delta=0,\quad\tau_{n+1}^\delta =
    \inf\set*{t\geq\tau_n^{\delta}}{\abs*{B_t-B_{\tau^\delta_n}}\geq\delta},
    \quad n\geq0.
  \end{displaymath}
  With this choice, the sequence $S_n^{\delta}=B_{\tau_n^\delta}$,
  $n\geq0$ forms a symmetric random walk 
  with step size $\delta$.
  
  The discrete \ito/ formula gives us that
  \begin{align*}
    F(S^\delta_{n})=F(0)
    &+\sum_{k=0}^{n-1}\frac{F(S^\delta_k+\delta)-F(S^\delta_k-\delta)}
      {2\delta}(S^\delta_{k+1}-S^{\delta}_k)\\
    &+\frac12\sum_{k=0}^{n-1}(F(S^\delta_k+\delta)+F(S^\delta_k-\delta)-2F(S^\delta_k)).
  \end{align*}
  It can be also written in the \ito/--Tanaka form
  \begin{align*}
    F(S^\delta_{n})=F(0)
    &+\sum_{k=0}^{n-1}\frac{F(S^\delta_k+\delta)-F(S^\delta_k-\delta)}{2\delta}
      (S^\delta_{k+1}-S^{\delta}_k)\\
    &+\frac12\sum_{r\in\Z}\zjel*{F((r+1)\delta)+F((r-1)\delta)-2F(r\delta)}
      \ell^\delta(r\delta,\tau^{\delta}_n),
  \end{align*}
  where $\ell^\delta(x,t)=\sum_{k:\tau^\delta_k<t} \I{S^\delta_k=x}$ denotes the
  number of visits of $S^\delta$ to the site $x$ before $t$.
  
  \bigskip

  Suppose now that $F(B)$ is a semimartingale with decomposition
  $F(B)=M+A$, where $M$ is a local martingale and $A$ is a process of
  finite variation. Denote $V$ the total variation process of
  $A$. Then $V$ has continuous sample paths taking finite
  values. Consider the following stopping times
  \begin{displaymath}
    \rho_K=\inf\set*{t\geq 0}{\max(V_t,\abs*{B_t})\geq K},
  \end{displaymath}
  and for a fixed $\delta$
  \begin{displaymath}
    \eta_K = \inf\set*{\tau^\delta_k}{\tau^\delta_k\geq\rho_K}.
  \end{displaymath}
  We later choose $K$ to be sufficiently large.
  From the relation
  \begin{displaymath}
    \E{F(S^\delta_{k+1})-F(S^{\delta}_k)\given\F_{\tau^\delta_k}}
    = \frac12\zjel*{F(S^\delta_k+\delta)+F(S^\delta_k-\delta)-2F(S^\delta_k)}  
  \end{displaymath}
  and that $\tau^\delta_k<\eta_K$ happens exactly when 
  $\tau^\delta_k<\rho_K$ we get
  \begin{align*}
    &\I{\tau^\delta_k<\eta_K}
      \frac12\zjel*{F(S^\delta_k+\delta)+F(S^\delta_k-\delta)-2F(S^\delta_k)}\\
    &=\I{\tau^\delta_k<\rho_K}\E{F(S^\delta_{k+1})-F(S^{\delta}_k)|\F_{\tau^\delta_k}}\\
    &=\I{\tau^\delta_k<\rho_K}
      \E{F(B_{\tau^\delta_{k+1}})-F(B_{\tau^\delta_{k+1}\wedge\rho_K})+
        F(B_{\tau^\delta_{k+1}\wedge\rho_K})-F(B_{\tau^\delta_{k}\wedge\rho_K})
        \given \F_{\tau^\delta_k}}\\
    &=\I{\tau^\delta_k<\rho_K}
      \E{F(B_{\tau^\delta_{k+1}})-F(B_{\tau^\delta_{k+1}\wedge\rho_K})\given \F_{\tau^\delta_k}}+
      \E{A_{\tau^{\delta}_{k+1}\wedge\rho_K}-A_{\tau^\delta_k\wedge\rho_K}\given \F_{\tau^\delta_k}}.
  \end{align*}
  Estimating the increment of $A$ with that of $V$ and using the boundedness of $F(B^{\eta_K})$
  this leads to
  \begin{multline*}
    \I{\tau^\delta_k<\eta_K}
    \frac12\abs*{F(S^\delta_k+\delta)+F(S^\delta_k-\delta)-2F(S^\delta_k)}\\
    \leq
    \E{V_{\tau^\delta_{k+1}\wedge\rho_{K}}-V_{\tau^{\delta}_{k}\wedge\rho_K}\given F_{\tau^\delta_k}}+
    2c(K,\delta)\P{\tau^\delta_{k+1} \geq \rho_K>\tau^\delta_k\given \F_{\tau^\delta_k}},
  \end{multline*}
  where
  \begin{displaymath}
    c(K,\delta) = \sup_{\abs*{x}<K+\delta}\abs*{F(x)}.
  \end{displaymath}
  Taking expectation of the sum we get for $\delta<1$ 
  that
  \begin{displaymath}
    \E{\sum_{r\in\Z}
      \frac12\abs*{F((r+1)\delta)+F((r-1)\delta)-2F(r\delta)}
      \ell^\delta(r\delta,\eta_K)}
    \leq K+2c(K,1).
  \end{displaymath}
  A simple calculation shows that
  \begin{displaymath}
    \E{L^{r\delta}_{\tau^\delta_{k+1}}-L^{r\delta}_{\tau^\delta_{k}}}
    =\delta\E{\I{S^{\delta}_{k}=r\delta}},
  \end{displaymath}
  hence
  \begin{displaymath}
    \E{L^{r\delta}_{\rho_K}}\leq
    \E{L^{r\delta}_{\eta_K}}
    = \delta\E{\ell^\delta(r\delta,\eta_K)}.
  \end{displaymath}
  From this we obtain that
  \begin{displaymath}
    \sum_{r\in\Z}
      \frac1{2\delta}\abs*{F((r+1)\delta)+F((r-1)\delta)-2F(r\delta)}
      \E{L^{r\delta}_{\rho_K}}
    \leq K+2c(K).
  \end{displaymath}
  
  Finally, let $I$ be a bounded interval.
  Since $\rho_K\to\infty$ as $K\to\infty$, there is a $K$ such
  \begin{displaymath}
    \inf_{x\in I} \E{L^x_{\rho_K}}>0.
  \end{displaymath}
  To see that this is really the case, note that $x\mapsto
  \E{L^x_{\rho_K\wedge1}}=\E{\abs*{B_{\rho_K\wedge 1}-x}-\abs*{x}}$  is
  a continuous function
  for all $K$ and as $K$ goes to infinity they tend to the everywhere
  positive function
  $x\mapsto\E{L^x_1}=\E{\abs*{B_{1}-x}-\abs*{x}}$ in a pointwise
  increasing manner,
  so by virtue of the Dini lemma the convergence is uniform
  on compacts. 

  Then
  \begin{displaymath}
    \mu_\delta(I)=\sum_{r:r\delta\in I}\frac1\delta
    \abs*{F((r+1)\delta)+F((r-1)\delta)-2F(r\delta)}
    \leq 2\frac{K+2c(K)}{\inf_{x\in I}\E{L^x_{\rho_K}}}<\infty.
  \end{displaymath}
  Since this upper estimate is independent of $\delta\in(0,1)$, we
  obtain that
  \begin{displaymath}
    \sup_{0<\delta<1}\mu_{\delta}(I)<\infty,
  \end{displaymath}
  which together with  Proposition \ref{prop:tot var} proves the necessity of our condition.
\end{proof}

Next, let $X$ be a continuous semimartingale, and suppose that $F$ is
locally Lip\-schitz.
The next theorem provides a necessary condition for a deterministic
function $F$ to be a semimartingale function for a continuous
semimartingale $X$.
\begin{theorem}\label{thm:X}
  Let $X$ be a continuous semimartingale in the filtration
  $\F$,
  denote the local time of $X$ by $(L^x_t)_{x\in\real,t\geq0}$ and put
  \begin{displaymath}
    H(X)=\set*{x\in\real}{\lim_{t\to\infty}\min(\E{L^x_t},{\E}\zjel*{L^{x-}_t})>0}.
  \end{displaymath}
  Suppose that $F$ is locally Lipschitz continuous.

  If $F(X)$ is a semimartingale in $\F$, then each point in $H(X)$ has a neighborhood $I$ 
  such that $F$ is the difference of   
  convex  functions on $I$.

  In particular, if $F$ is continuously differentiable, then
  each point in $H(X)$ has a neighborhood $I$ 
  such that the total variation of $F'$ is finite 
  on $I$.
\end{theorem}

\begin{proof}
  We are going to use a similar argument as in the case of the
  Brownian motion. Suppose that $X=M^X+A^X$, where $M^X$ is a continuous local
  martingale and $A^X$ is a process of finite variation. Denote $V^X$
  the total variation process of $A^X$. As $F (X)$ is a
  semimartingale,
  \begin{displaymath}
    F(X_t)=F(X_0)+M_t+A_t,
  \end{displaymath}
  where $M$ is a local martingale in the filtration $\F$ with
  continuous sample paths starting from zero and $A$ is a 
  process of finite variation also starting from zero.
  Denote the total variation process of $A$ by $V$. Let us consider the
  following stopping time
  \begin{displaymath}
    \rho_K=K\wedge\inf
    \set*{t\geq0}
    {\max(V_t,V_t^X,\abs*{X_t}
      ,\br*{M}_t
      )\geq K},
  \end{displaymath}
  where $K$ will be chosen sufficiently large. With this choice
  the local martingale part $(M^X)^{\rho_K}$ of the stopped process $X^{\rho_K}$ is a true
  martingale and $V^{\rho_K}$ is integrable.

  Since we want to apply the discrete Ito formula to the embedded
  random walk, we modify $X$ after $\rho_K$.  
  By enlarging the probability space we can assume that
  there is a Brownian motion $B$ independent of $\F_\infty$. Let
  \begin{align*}
    \tX_t&= X_{t\wedge\rho_K}+B_{t}-B_{t\wedge\rho_K},\\
    \tF_t&=\F_t\vee \F^B_t.
  \end{align*}
  Shortly, we can assume 
  that for each $\delta>0$ the stopping times
  \begin{align*}
    \tau^\delta_0
    &=\inf\set*{t\geq0}{\tX_t\in \delta\Z},\\
    \tau^\delta_{k+1}
    &=\inf\set*{t\geq0}{\abs*{\tX_t-\tX_{\tau^\delta_k}}=\delta},
      \quad k\geq0,\\
    \eta_K
    &= \inf\set*{\tau^\delta_k}{\tau^\delta_k\geq\rho_K},
  \end{align*}
  are all almost surely finite. As in the case of the Brownian
  motion $S^\delta_k=\tX_{\tau^\delta_k}$ is a random walk on the
  lattice $\delta\Z$, although not necessarily symmetric. As in the proof of
  Theorem \ref{thm:BM}
  \begin{align*}
    &\I{\tau^\delta_k<\eta_K}
      \frac12\zjel*{F(S^\delta_k+\delta)+F(S^\delta_k-\delta)-2F(S^\delta_k)}\\
    &=\I{\tau^\delta_k<\rho_K}
      \E{F(S^\delta_{k+1})-F(S^{\delta}_k)-
      \frac{F(S^\delta_k+\delta)-F(S^\delta_k-\delta)}{2\delta}
      (S^\delta_{k+1}-S^\delta_k)\given \F_{\tau^\delta_k}}.
  \end{align*}
  This identity follows from the fact that the left hand side is
  $\F_{\tau^\delta_{k}}$ measurable. In this case, however,
  \begin{displaymath}
    \I{\tau^\delta_k<\rho_K}\E{S^\delta_{k+1}-S^\delta_k\given \F_{\tau^\delta_k}}=
    \E{A^X_{\tau^\delta_{k+1}\wedge\rho_K}-A^X_{\tau^\delta_k\wedge\rho_K}\given \F_{\tau^\delta_k}},
  \end{displaymath}
  and
  \begin{multline*}
    \I{\tau^\delta_k<\rho_K}
    \E{F(S^\delta_{k+1})-F(S^\delta_{k})\given \F_{\tau^{\delta}_k}}\\
    =
    \I{\tau^\delta_k<\rho_K}
    \E{F(\tX_{\tau^\delta_{k+1}})-F(\tX_{\tau^\delta_{k+1}\wedge\rho_K})\given \F_{\tau^{\delta}_k}} 
    +
    \E{F(\tX_{\tau^\delta_{k+1}\wedge\rho_K})-F(\tX_{\tau^\delta_{k}\wedge\rho_K})\given \F_{\tau^{\delta}_k}}.
  \end{multline*}
  Here
  \begin{align*}
    \E{F(\tX_{\tau^\delta_{k+1}\wedge\rho_K})-
    F(\tX_{\tau^\delta_{k}\wedge\rho_K})\given \F_{\tau^{\delta}_k}}=
    \E{
    A_{\tau^\delta_{k+1}\wedge\rho_K}-A_{\tau^\delta_{k}\wedge\rho_K}
    \given \F_{\tau^{\delta}_k}}.
  \end{align*}
  From these we get that
  \begin{multline*}
    \E{\I{\tau^\delta_k<\rho_K}
      \frac12\abs*{F(S^\delta_k+\delta)+F(S^\delta_k-\delta)-2F(S^\delta_k)}}\\
    \leq
    c(K,\delta)
    \E{V^X_{\tau^{\delta}_{k+1}\wedge\rho_K}-V^X_{\tau^{\delta}_{k}\wedge\rho_K}}
      +\E{V_{\tau^{\delta}_{k+1}\wedge\rho_K}-V_{\tau^{\delta}_{k}\wedge\rho_K}}
      +2C(K,\delta)\P{\tau^\delta_k<\rho_K\leq\tau^{\delta}_{k+1}}
  \end{multline*}
  where
  \begin{displaymath}
    c(K,\delta)=\sup\set*{\abs*{\frac{F(y)-F(x)}{y-x}}}{\abs*{x},\abs*{y}\leq K+\delta,\, x\neq y}, 
    \quad
    C(K,\delta)=\max_{\abs*{x}\leq K+\delta}\abs*{F(x)}.
  \end{displaymath}
  After summation
  \begin{multline}
    \label{eq:mu delta}
    \sum_{r\in\Z} 
    \frac12\abs*{F((r+1)\delta)+F((r-1)\delta)-2F(rk)}\E{\ell^\delta(r\delta,\rho_K)}\\
    \leq
    K(1+2c(K,\delta))+2C(K,\delta)<\infty. 
  \end{multline}
  To finish the proof we use Lemma \ref{lem:ell} below, which
  for each $x\in H(X)$ provides us a
  $K>0$, a non-empty open  interval $I$ and $\delta_0>0$, such that
  \begin{equation*}
    \inf_{\delta\in (0,\delta_0)}
    \inf_{r:r\delta\in I}\delta\E{\ell^{\delta} (r\delta,\rho_K)}>0.
  \end{equation*}
  Then a rearrangement of \eqref{eq:mu delta} gives 
  \begin{displaymath}
    \limsup_{\delta\to0}\mu_{\delta}(I) <\infty,
  \end{displaymath}
  which by Proposition \ref{prop:tot var} proves the claim.
\end{proof}

\section{Random semimartingale function}\label{sec:rand F}

Throughout this section $(S,\cS)$ will be a measurable space and $\xi$
will be 
a random variable taking values in $S$. As it can cause no confusion,
for the sake of brevity, we use ``for almost all $z\in S$'' instead 
of the correct but longer form  ``for $\P\circ\xi^{-1}$ almost all $z\in S$''. 
We will consider random
variables (and processes) of the form $\tX (\xi)$
where $\tX:\tOmega=\Omega\times S \to\real$. The goal here is to show
a version  of Theorem \ref{thm:X} when the function $F$ is also
random. For the precise formulation we replace the single variable
function  $F$  with a parametric version $F:\real\times S\to
\real$, and consider the process $(F (X_t,\xi))_{t\geq 0}$. This has
the form $(\tX_t (\xi))_{t\geq 0}$ with $\tX_t(z)=F(X_t,z)$. 
Suppose now that $\tX (\xi)$ is a
semimartingale in the enlarged filtration $\F^\xi$, where
$\F^\xi_t=\F_t\vee\sigma (\xi)$ and $\xi$ is independent of $\F_\infty$. 
Note that the right continuity of $\F$ is inherited to  $\F^\xi$ by the 
independence of $\F_\infty$ and $\xi$. 
We show in Theorem \ref{thm:rand F}
below that in this case $\tX (z)$ is also a semimartingale in $\F$ for
almost all $z\in S$.  So, if $(F (X_t,\xi))_{t\geq0}$ is a
semimartingale in $\F^\xi$, then $F (X_t,z)$ is a semimartingale in
$\F$ for almost all $z$ and  Theorem \ref{thm:X} applies to $x\mapsto F
(x,z)$ for almost all $z$.

\begin{theorem}\label{thm:rand F}
  Let $\F$ be a filtration such that $\xi$ is independent of
  $\F_\infty$.
  
  Suppose that $\tX:[0,\infty)\times \tOmega\to\real$ is such
  that $\tX (\xi)$ is a continuous semimartingale in $\F^\xi$. Then
  $\tX (z)$ is a semimartingale in $\F$ for almost all $z\in S$.
\end{theorem}

\begin{corollary}\label{cor:rand F}
  Let $\F$ be a filtration such that $\xi$ is independent of
  $\F_\infty$. Suppose that $X$ is a continuous semimartingale in $\F$,
  $F:\real\times S\to\real$ is locally Lipschitz in the first variable 
  and $(F (X_t,\xi))_{t\geq0}$ is a semimartingale in $\F^\xi$. 
  
  Then each point in $H(X)$ has a neighborhood $I(z)$ for almost all $z\in S$ 
  such that $x\mapsto F (x,z)$ is the difference of   
  convex  functions on $I(z)$.

  In particular, if $F$ is continuously differentiable in its first variable 
  for almost all $z$, then the following holds for almost all $z$:
  each point in $H(X)$ has a neighborhood $I(z)$ 
  such that the total variation of $x\mapsto\partial_x F(x,z)$ is finite 
  on $I(z)$.
\end{corollary}

We start with preliminary lemmas. In the next lemmas for a $\sigma$ algebra $\cA$  the family of $\cA$ 
measurable functions is denoted by $m(\cA)$. 

\begin{lemma}\label{lem:par X rv}
  Let $\cA$ be a $\sigma$-algebra of events and $X\in m(\cA\vee\sigma (\xi))$.
  Then $X=\tX (\xi)$, where
  $\tX:\tOmega\to\real$ and $\tX$ is $\cA\times\cS$ measurable. 
\end{lemma}
\begin{proof}
  This is  a monotone class argument. Let $\cH$ be the collection of
  random variables having the desired representation property, that is
  \begin{displaymath}
    \cH=\set*{Y\in m(\cA\vee\sigma (\xi))}
             {\exists \tY\in m (\cA\times\cS),\,Y=\tY (\xi)}.
  \end{displaymath}
  $\cH$ is obviously a linear space. Although $\tY$ is not unique one
  can find $\tY$  for $Y\in\cH$, such that $\sup \abs{Y}=\sup\abs{\tY}$ and therefore
  $\cH$  is closed under uniform convergence.
  
  Then let
  \begin{displaymath}
    \cD=\set*{A\in\cA\vee\sigma (\xi)}{\I[A]\in\cH}.
  \end{displaymath}
  $\cD$ is clearly a $\lambda$-system, and it is also obvious that
  $\cD$ extends the $\pi$-system
  \begin{displaymath}
    \cC=\set*{A\cap B}{A\in\cA,\,B\in\sigma (\xi)},
  \end{displaymath}
  which generates $\cA\vee\sigma (\xi)$. So $\cD=\cA\vee\sigma (\xi)$
  and $\cH$ contains all bounded $\cA\vee\sigma (\xi)$ measurable
  random variables. Now if $X\in m (\cA\vee\sigma (\xi))$, then $Y=\arctg
  (X)\in \cH$  and $\tX=\I{\abs{\tY}<\pi/2}\tg (\tY)$ shows that $X\in
  \cH$. 
\end{proof}

\begin{lemma}\label{lem:par X proc}
  Let  $\F$ be a filtration.  
  Suppose that $X$ is a \cadlag/ process adapted to the filtration
  $(\F^\xi_t)_{t\geq0}$.
  Then there is a parametric process
  $\tX:[0,\infty)\times \tOmega\to\real$, such that
  \begin{enumerate}
  \item $X$ and $\tX(\xi)$ are
    indistinguishable,
  \item $\tX (z)$ is adapted to $\F$ for all $z\in S$. 
  \end{enumerate}
  When $\F_\infty$ is independent of $\xi$ then 
  the sample path properties are inherited by $\tX$, that is
  \begin{enumerate}
  \setcounter{enumi}{2}
  \item $\tX (z)$ is \cadlag/ for  almost all $z\in S$.
  \item If $X$ have continuous sample paths, the same is true for $\tX(z)$
  for almost all $z\in S$
  \item If $X$ is of finite variation, then the same is true for $\tX(z)$ 
  for almost all $z\in S$. 
  \end{enumerate}
\end{lemma}
\begin{proof}
  Via Lemma \ref{lem:par X rv} there are $\tX_q\in m(\F_q \times \cS)$
  for $q\in[0,\infty)\cap\Q$ such that $\tX_q(\xi)=X_q$. Then
  \begin{displaymath}
    \tX_t(\omega,z)=\liminf_{q\searrow t}\tX_q(\omega,z),\quad t\geq0
  \end{displaymath}
  defines for each $t\geq0$ an $\F_t\times\cS$ measurable
  function.
  Here we used the right continuity of $\F$, which is
  tacitly assumed. In other words, $\tX(z)$ is adapted to $\F$ for all
  $z\in S$.
  
  From the definition of $\tX$ we also
  have that $\tX_t(\xi)=X_t$ for all $t\geq0$ whenever the sample path of
  $X$ is
  \cadlag/, which happens almost surely, hence
  $\tX(\xi)$ and $X$ are indistinguishable.

  Suppose now that $\F_\infty$ and $\xi$ are independent.
  The last part of the claim is the application of the Fubini theorem. 
  Let $A\subset \Omega$ be the collection of outcomes $\omega\in \Omega$ 
  such that $(X_q(\omega))_{q\in \Q\cap[0,\infty)}$ has right limits 
  at each $t\geq 0$ and left limits at each $t>0$. 
  Then, using for example the upcrossing number, 
  $A\in\sigma(\set{X_q}{q\in\Q\cap[0,\infty)})$. Similarly, 
  let $g:S\times\Omega\to\smallset{0,1}$ the indicator that 
  $(\tX_q(z))_{q\in\Q\cap[0,\infty)}$ has right and left limits at each $t$. 
  The function $g$ satisfies $g(\xi)=\I[A]$ and $g\in m(\cS\times\cA)$. 
  As by assumption $\P{A}=1$, so we have that $\E{g(z,\cdot)}=1$ for almost all $z\in S$. 
  That is, for almost all $z$, $\tX(z)$, which is the right continuous extension 
  of $(\tX_q(z))_{q\in\Q\cap[0,\infty)}$, is \cadlag/.

  The last two parts of the claim go along the same line, so they are left to the reader.
\end{proof}

\begin{lemma}\label{lem:par X cond}
  Let $\G\subset\F$ be two $\sigma$-algebras.
  Suppose that $\xi$ is independent of $\F$.
  Using Lemma \ref{lem:par X rv} write $X\in L^1(\F\vee\sigma(\xi))$ and
  $Y=\E{X\given \G\vee\sigma(\xi)}$ as 
  $X=\tX (\xi)$ and $Y=\tY (\xi)$.
  Then
  $\tY(z)=\E{\tX(z)\given \G}$
  almost surely for almost all $z\in S$.
\end{lemma}
\begin{proof}
  $X$ is integrable, so by the independence of $\xi$ and $\F$ we also
  have that $\tX (z)$ is integrable for 
  almost all $z \in S$. 
  By the monotone class argument $\E{\tX(z)|\G}$ has a version which is
  $\G\times\cS$ measurable, that is there is $\tZ\in m(\G\times\cS)$
  such that $\E{\tX(z)|\G}=\tZ(z)$ almost surely for 
  almost all $z \in S$.
  
  Let
  \begin{displaymath}
    A=\event{\tZ (\xi)>\tY(\xi)}\quad\text{and}\quad
    A(z)=\event{\tZ (z)>\tY(z)}.
  \end{displaymath}
  Then $A\in\G\vee\sigma(\xi)$ and $\I[A]=\I[A(z)]|_{z=\xi}$, so by
  the independence of $\xi$ and $\F$
  \begin{align*}
    \E{(\tX(\xi)-\tY(\xi))\I[A]\given \xi}
    &=\left.\E{(\tX(z)-\tY(z))\I[A(z)]}\right|_{z=\xi}\\
    &=\left.\E{(\tZ(z)-\tY(z))\I[A(z)]}\right|_{z=\xi}.
  \end{align*}
  In the last step we used that $\E{\tX (z)\given \G}=\tZ (z)$.
  
  Then as $Y=\E{X\given \G\vee\sigma(\xi)}$ we have that
  \begin{displaymath}
    0=\E{(X-Y)\I[A]}=\E{g(\xi)},
    \quad\text{where}\quad
    g(z)=\E{(\tZ(z)-\tY(z))\I[A(z)]}.
  \end{displaymath}
  Since $g$ is non negative by the choice of $A (z)$, we have that $g$
  is zero for almost all  $z$ and 
  $\tZ(z)\leq \tY(z)$ almost
  surely for almost all $z$. The other direction is obtained
  similarly.
  
  So for almost all $z$ the 
  $\G$ measurable random variables
  $\tY(z)$, $\tZ(z)$ and $\E{\tX(z)\given \G}$ are almost surely equal,
  which proves the claim.
\end{proof}

\begin{lemma}\label{lem:par X martingale}
  Let  $\F$ be a filtration.
  Suppose that $\xi$ is independent of $\F_\infty$ and
  $X$ is a \cadlag/ martingale in the filtration $(\F^\xi_t)_{t\geq0}$.

  Then there is a parametric process
  $\tX:[0,\infty)\times\tOmega\to\real$ such that $\tX(\xi)$
  and $X$ are indistinguishable and $\tX(z)$ is a martingale in $\F$
  for almost all $z$.
\end{lemma}
\begin{proof}
  Let $\tX$ be the parametric process obtained from the application of
  Lemma \ref{lem:par X proc}. By Lemma \ref{lem:par X cond} $(\tX_t
  (z))_{t\in[0,\infty)\cap \Q}$ is martingale in
  $(\F_t)_{t\in[0,\infty)\cap \Q}$ for $z\in H$ where $\P{\xi\in
    H}=1$. As we have seen in Lemma \ref{lem:par X proc} $\tX(z)$
  is a \cadlag/ process for almost all $z$, hence for $z\in H'$ with
  $\P{\xi\in H'}=1$.

  As the filtration $\F$ is right continuous we have for $z\in H\cap
  H'$ the \cadlag/ process $\tX(z)$ is a martingale in $\F$.
\end{proof}

The localized version of Lemma \ref{lem:par X martingale} follows
from the next claim.
\begin{lemma}\label{lem:par tau}
  Suppose that $\tau$ is a stopping time in $\F^\xi$. Then there is
  $\ttau:\tOmega\to[0,\infty]$ such that $\ttau (\xi)=\tau$
  and $\ttau (z)$ is a stopping time in $\F$ for
  all $z\in S$. 
\end{lemma}
\begin{proof}
  $\event{\tau< t}\in \F^\xi_t$, so using Lemma \ref{lem:par X rv}, 
  $\I{\tau< t}=\tY_t(\xi)$ where $\tY_t\in m(\F_t\times\cS)$.
  If we define
  \begin{displaymath}
    \ttau(z)=\inf\set*{q\in\Q\cap(0,\infty)}{\tY_q(z)\neq 0},\quad z\in S,
  \end{displaymath}
  then, as for fixed $z$ the random variable $\tY_q(z)\in m(\F_q)$, it follows 
  that $\event{\ttau(z)<t}\in\F_t$ for all $t\geq0$, which by the right 
  continuity of $\F$ implies that $\ttau(z)$ is a stopping time for all $z\in S$.
  It is an easy exercise that 
  \begin{displaymath}
    \event{\ttau(\xi)<t}=\cup_{q<t}\event{\tY_q(\xi)\neq 0}=\cup_{q<t}\event{\tau<q}=\event{\tau<t},
  \end{displaymath}
  and $\ttau(\xi)=\tau$ follows.
\end{proof}

Now we turn to the 
\begin{proof}[Proof of Theorem \ref{thm:rand F}]
  $\tX(\xi)$ is a semimartingale in $\F^\xi$, which means that $\tX(\xi)=M+A$, 
  where $M$ is local martingale and $A$ is a process of finite variation.
  Then $M=\tM(\xi)$ and $A=\tA(\xi)$ by Lemma \ref{lem:par X proc}, moreover
  both $\tM(z)$ and $\tA(z)$ inherit the path properties of 
  $M$ and $A$, respectively, for almost all $z\in S$.
  
  Combining Lemma \ref{lem:par X martingale} 
  and Lemma \ref{lem:par tau} we get that $\tM(z)$ is a local martingale 
  in $\F$ for almost all $z\in S$. 

  Finally, using the Fubini theorem we get that $\tX(z)$ and 
  $\tM(z)+\tA(z)$ are indistinguishable for almost all $z\in S$.
\end{proof}

\section{Examples}\label{sec:examples}

In the following few examples $S$ is either $C(\real)$ or $C([0,s])$ for some $s>0$, 
that is, the space of continuous functions. With a suitable metric $S$ 
is complete separable metric space and its Borel $\sigma$-algebra $\cS$ is the same 
as a the smallest $\sigma$-algebra making all coordinate mapping measurable. 
We consider two type of examples here 
\begin{displaymath}
  F(x,w)=\int_0^x w_sds \quad\text{or}\quad F(x,w)=\int_0^s \I{w(u)<x}du.
\end{displaymath}
In the first version $x\mapsto F(x,\xi)$ is continuously differentiable 
regardless of the choice of $\xi$,
in the second version it is continuously differentiable almost surely,
for example, when $\xi$ is a Brownian motion.

Suppose that $\xi$ is a random variable taking values in $S$ and 
independent of the Brownian motion $B$. Denote $\F$ the natural filtration of $B$.
Then since $H(B)=\real$, for the process $(F(B_t,\xi))_{t\geq 0}$ 
to be a semimartingale in $\F^\xi$ it is needed that 
$x\mapsto F(x,\xi)$ is the difference of two convex functions almost surely. 
When $F(x,\xi)$ is $C^1$  in $x$, it 
simply requires that the sample paths of $x\mapsto \partial_x F(x,\xi)$ have to have finite 
total variation on compact intervals.

\begin{proposition} \label{prop_bm}
  Let $B$ and $\xi$ be two independent standard Brownian motions and 
  let $\F = (\F_t)_{t \geq 0}$ be the natural filtration of $B$ and $\F^\xi$ is 
  the initial enlargement of $\F$ with $\sigma(\xi)$, that is
  $\F^\xi_t = \F_t\vee\sigma(\xi)$, $t\geq 0$. 
  Let $F$ be the following random function:
  \begin{displaymath}
    F(x,\xi) = 
    \begin{cases}
      \int_{0}^x \xi_{y} d y, \quad x \geq 0, \\
      0, \quad x < 0.
    \end{cases}
  \end{displaymath}
  Then $(F(B_t,\xi))_{t\geq0}$ is not an $\F^\xi$--semimartingale.
\end{proposition}
\begin{proof}
  In this example $\partial_x F(x,\xi)=\xi_x$, where $\xi$ is a Brownian motion.
  Since the total variation of $\xi$ is almost surely infinite on any non-empty sub-interval 
  of the positive half line, $x\mapsto F(x,\xi)$ is not a difference of convex 
  functions almost surely and 
  $(F(B_t,\xi)$ is not a semimartingale. 
\end{proof}

\begin{proposition} \label{prop_fbm} 
  Let $\xi=B_H$ be a fractional Brownian motion with Hurst index $H
  \in (0,1)$ and $B$ be a Brownian motion, independent of
  $B_H$. Let $\F = (\F_t)_{t \geq 0}$ be the filtration
  generated by $B$ and $\F^\xi$ be the initial enlargement of $\F$ with $\sigma(B^H)$, that is $\F^\xi_t = \F_t
  \vee \sigma(B^H)$, $t \geq 0$. Let $F$ be the following function:
  \begin{displaymath}
    F(x,\xi) = 
    \begin{cases}
      \int_0^x \xi_y dy,& x \geq 0, \\
      0,& x < 0.
    \end{cases}
  \end{displaymath}
  Then $(F(B_t,\xi))_{t\geq 0}$  is not an $\F^\xi$--semimartingale.
\end{proposition}
\begin{proof} 
  It is well know that the total variation of the fractional 
  Brownian motion  is almost surely infinite on any nonempty sub-interval of the positive half 
  line, see for example \cite{pratelli} and the references therein, 
  and the statement follows the same way as in the case of Proposition \ref{prop_bm}
\end{proof}

As an application of Theorem \ref{thm:rand F} we can obtain a
preparatory example for the  process $A(t,B_t)$ investigated in \cite{MR1118450}.
Here $A(t,x)$ is the amount of time spent by the Brownian motion $B$ below the 
level $x$ up to time $t$. 
First, we treat the case where there are two independent Brownian motions involved.
\begin{proposition}\label{prop:prep Rogers}
  Let $B$  and $\xi$ be two independent standard Brownian motions. 
  For a fixed $s>0$ consider the following function
  \begin{displaymath}
    F(x, \xi) = \int_0^s \I{\xi_u \leq x}d u. 
  \end{displaymath}
  Let $\F$ be the natural filtration of $B$ 
  and $\F^\xi$ be the initial enlargement of $\F$ with $\sigma(\xi)$, 
  that is $\F^\xi_t = \F_t \vee\sigma(\xi), t \geq 0$.
  
  Then $(F(B_t, \xi))_{t \geq 0}$ is not an $\F^\xi$--semimartingale.
\end{proposition}
\begin{proof}
  By the occupation time formula
  \begin{displaymath}
    F(x,\xi) = \int_{-\infty}^x Z_ydy,
  \end{displaymath}
  where $Z_y=L^y_s(\xi)$ denotes the local time profile of $\xi$ at time $s$.
  For a compact interval $[a,b]$ the quadratic variation of $Z$ 
  (see \cite[Theorem (1.12)]{revuz-yor}) exists and is given by the following formula
  \begin{displaymath}
    \br{Z}_b - \br{Z}_a =\int_a^b 4 Z_y dy.
  \end{displaymath}
  Since $P(L_s^0(\xi)>0)=1$ 
  and $Z_y = L_s^y(\xi)$ is continuous in $y$, we can conclude that 
  $\P{\int_a^b Z_ydy > 0}=1$ also holds for any $a<0<b$. 
  From this it easily follows that on an almost sure event 
  $x\mapsto Z_x$ has infinite total variation on any non-empty interval 
  containing zero. So almost surely $F$ is not a semimartingale function for $B$ and
  $(F(B_t,\xi))_{t \geq 0}$
  can not be a semimartingale in $\F^\xi$ by Theorem \ref{thm:rand F}.  
\end{proof}

Now we turn to the process $A(t,B_t)$ of \cite{MR1118450}.  
\begin{proposition}\label{prop:Rogers}
Let $B$ be a standard Brownian motion and denote $\F$  its natural filtration. 
For $x\in\real$ put 
\begin{displaymath}
  A(t,x) = \int_0^t \I{B_u \leq x}du. 
\end{displaymath}
Then $A(t,B_t)$ is not semimartingale in $\F$.
\end{proposition}
The result in \cite{MR1118450} actually states more than just the lack of the 
semimartingale property, namely, they define
\begin{displaymath}
  X_t=A(t,B_t)-\int_0^t L_s^{B_s} dB_s
\end{displaymath}
and show that the $p$ order variation of $X$ on $[0,1]$ is zero when $p>4/3$ and infinite 
when $p<4/3$.
\begin{proof}
  Fix $s > 0$ and for $t\geq s$ write $A$ as 
  \begin{displaymath}
    A(t,x) = \int_0^t \I{B_u \leq x}du = 
    \int_0^s \I{B_u \leq x}du + \int_s^{t} \I{B_{u}-B_s \leq x-B_s}du.
  \end{displaymath}
  Denote by $\beta_t = B_{s+t}-B_{s}$, $t\geq 0$. Then $\zjel{\beta_u}_{u \geq 0}$ is 
  a Brownian motion independent
  of $\F_s$.

  With 
  \begin{displaymath}
    \alpha(t,x) = \int_{0}^t \I{\beta_u \leq x}du 
  \end{displaymath}
  we can write for $t\geq0$ 
  \begin{equation}\label{eq:A s+t}
    A(s+t,x) = A(s,x) + \alpha(t, x-B_s), \quad 
    A(s+t,B_{s+t}) = A(s,B_{s+t}) + \alpha(t,\beta_{t}).
  \end{equation}
  
  Suppose now, on the contrary to the claim that $(A(t,B_t))_{t\geq 0}$ is a semimartingale
  in the natural filtration of $B$, that is, in $\F$. 
  Then $\zjel*{\alpha(t, \beta_{t})}_{t \geq 0}$ is a semimartingale in the 
  filtration of $\beta$. As $(\G_t=\F_{s+t})_{t\geq 0}$ is an initial enlargement 
  of the natural filtration of $\beta$ with an independent $\sigma$-algebra $\F_s$, 
  the process $\zjel*{\alpha(t, \beta_{t})}_{t \geq 0}$ is also a semimartingale in 
  $(\G_t=\F_{s+t})_{t\geq 0}$.  Since  $(A(s+t,B_{s+t}))_{t\geq 0}$ is obviously 
  a semimartingale in $(\G_t)_{t\geq0}$, we get  from the second part of
  \eqref{eq:A s+t} that $(A(s,B_{s+t}))_{t\geq 0}$ also needs 
  to be a $(\G_t)_{t\geq 0}$ semimartingale.
  But
  \begin{displaymath}
    A(s,B_{s+t}) = \int_0^s \I{B_u \leq B_{t+s}} du =
    \int_0^s \I{B_u - B_s \leq B_{t+s} - B_s}du. 
  \end{displaymath}
  Note that $\xi_u = B_{s-u} - B_s$, $u\in[0,s]$, which is the time reversal of 
  $(B_u)|_{u\in[0,s]}$, is a Brownian motion (on $[0,s]$),
  independent of $\zjel{\beta_{t}}_{t\geq 0}$. 
  With
  \begin{displaymath}
    F(x,\xi) = \int_0^s \I{\xi_u \leq x}du=\int_0^s \I{B_u-B_s \leq x}du, 
  \end{displaymath}
  we have 
  \begin{displaymath}
    A(s,B_{s+t}) = F(\beta_{t},\xi).
  \end{displaymath} 
  In Proposition \ref{prop:prep Rogers} we already proved that 
  $(F(\beta_t,\xi))_{t\geq 0}$ is not a semimartingale in the initial enlargement 
  of the natural filtration of $\beta$ with $\xi$ which is the same as 
  $\G$. On the other hand, the assumption that $(A(t,B_t))_{t\geq0}$ is a 
  semimartingale in $\F$  would lead to the conclusion that 
  $(A(s,B_{s+t})_{t\geq0}$ is semimartingale in $\G$. This
  contradiction proves the claim.
\end{proof}

\section{The median process is not a semimartingale}\label{sec:m not semi}

Consider the equation
\begin{equation}\label{eq:D (x)}
  d D_t (x)=\sigma (D_t(x))dB_t, \quad D_0(x)=x, 
\end{equation}
where $\sigma (x)=x\wedge (1-x)$. The aim of this section is to show that 
$(D_t^{-1}(\alpha))_{t\geq0}$ is not a semimartingale for $\alpha\in(0,1)$, 
in particular, the conditional median $m_t=D_t^{-1}(1/2)$ 
lacks the semimartingale property. It is useful to collect some of the properties of $D$.
It was shown in \cite{drift:2009} that $x\mapsto D_t(x)$ is differentiable for all $t$ 
almost surely and the space derivative is given by the stochastic exponential of
\begin{displaymath}
    \int_0^t \sigma'(D_s(x)) dB_s,
\end{displaymath}
provided that $\sigma$ is Lipshitz continuous. 
In this case $\sigma'$ exists Lebesgue almost everywhere and is a bounded function.
Note that, it also follows, for example from the positivity of $D'_t$, 
that $x\mapsto D_t(x)$ is stricly increasing for all $t$ on 
an almost sure event.
We prove below that $(x,t)\mapsto D'_t(x)$ is continuous whenever $\sigma$ is 
Lipshitz and $\sigma'$ is monotone, for example, $\sigma(x)=x\wedge(1-x)$ has 
this property. In other words, $D'$ is continuous almost surely if $\sigma$ is a convex 
(or concave) Lipschitz function, that is $\sigma(x)=\int_0^x \sigma'$ for some monotone 
bounded function denoted by $\sigma'$. The proof of the continuity in 
Lemma \ref{lem:cont of D'} uses the Lamperti representation of $D$, and at this step 
the monotonocity of $\sigma'$ is important. On the other hand, all other parts of our 
argument only use that $\sigma'$ is of bounded variation on compact intervals, that is 
$\sigma$ is a difference of two convex functions. Even this assumption seems 
to be too strict. In the proof below 
we compute the quadratic variation of $D'_t$ with respect to the space variable $x$ and 
use that the total variation is infinite if the quadratic variation is non-zero.
Finding another way of proving the unboundedness of the total variation would relax 
the condition imposed on $\sigma$. 

\begin{theorem}\label{thm:D(x)}
  Suppose that $\sigma$ is a convex (or concave) Lipschitz function, 
  $B$ is a Brownian motion 
  and $(D,B)$ satisfies 
  \eqref{eq:D (x)}. Let $\alpha\in\real$ such that 
  $\sigma(\alpha)\neq0$ 
  and denote $C$ the connected component of 
  $\real\setminus\event{\sigma=0}$
  containing $\alpha$. 
    
  If $\sigma'$ is not continuous on $C$,
  then $q_t=D_t^{-1}(\alpha)$ is not a semimartingale in the natural filtration of $B$.
\end{theorem}

Before the proof we discuss the case $\sigma(x)=x\wedge(1-x)$. Then 
\begin{displaymath}
  D_t(x)=
  \begin{cases}
  x\exp{B_t-\frac12t},& x\leq 0,\\
  1-(1-x)\exp{-B_t-\frac12t},&x\geq 1.\\  
  \end{cases}
\end{displaymath}
So for $\alpha\notin (0,1)$ the process $q_t(\alpha)=D_t^{-1}(\alpha)$ can be explicitly 
given 
\begin{displaymath}
  q_t(\alpha)=
  \begin{cases}
  \alpha\exp{-B_t+\frac12t},& \alpha\leq 0,\\
  1-(1-\alpha)\exp{B_t+\frac12t},&\alpha\geq 1,\\  
  \end{cases}
\end{displaymath}
and from this explicit form we see that it is a semimartingale. 
On the other hand, $C=(0,1)$ contains $1/2$, the point of 
discontinuity for $\sigma'$ and for $\alpha\in(0,1)$ the process 
$(q_t(\alpha))_{t\geq0}$ is not a 
semimartingale by the theorem.

\begin{proof} We suppose on the contrary to the claim that $q$ is a semimartingale
  and show that this leads to a contradiction.

  Let $s>0$. Then, as it was already indicated in the introduction,
  we can decompose $(q_{s+t}=D_{s+t}^{-1}(\alpha))_{t\geq0}$ as
  \begin{equation}\label{eq:D_s D_{s,t}}
    D_{s+t}^{-1}(\alpha) = D_s^{-1}(D_{s,t}^{-1}(\alpha)),
  \end{equation}
  where  $D_{s,t}(x)$ is the solution of \eqref{eq:D (x)} with $B$ replaced by 
  $B_{s,t}=B_{s+t}-B_{s}$ a ``Brownian motion that starts evolving at time $s$''. 
  Then $D_{s+t}(x)=D_{s,t}(D_s(x))$ and this is the rationale behind \eqref{eq:D_s D_{s,t}}.

  As by our assumption $q=D^{-1}(\alpha)$ is a semimartingale in the filtration of the 
  driving Brownian motion, so is $(q_{s,t}=D_{s,t}^{-1}(\alpha))_{t\geq0}$ 
  in the filtration of $(B_{s,t})_{t\geq0}$ denoted by
  $(\F_{s,t})_{t\geq0}$, where $\F_{s,t}=\sigma(\set*{B_{s,u}}{u\leq t})$. 
  Denote by $\F$ the filtration of $B$, that is 
  $\F_t=\sigma(\set*{B_u}{u\leq t})$. Then $(q_{s,t})_{t\geq0}$ remains a 
  semimartingale in the filtration $(\F_{s+t})_{t\geq 0}$ as 
  $\F_{s+t}=\F_s\vee\F_{s,t}$, that is, $(\F_{s+t})_{t\geq 0}$ is an initial 
  enlargement of $(\F_{s,t})_{t\geq0}$ with an independent $\sigma$-algebra $\F_s$.
  
  On the other hand, as  $q$ is supposed to be a semimartingale in $\F$, the process 
  $(q_{s+t})_{t\geq0}$ is 
  a semimartingale in $(\F_{s+t})_{t\geq0}$. By \eqref{eq:D_s D_{s,t}}
  $q_{s+t}=D_{s}^{-1}(q_{s,t})$ and the random function 
  $D^{-1}_{s}$ is a semimartingale function for $(q_{s,t})_{t\geq 0}$.

  Our necessary condition for $(D_s^{-1}(q_{s,t}))_{t\geq0}$ being a semimartingale 
  involves the set $H((q_{s,t})_{t\geq0})$ defined in Theorem \ref{thm:X}.
  Below we show in Corollary \ref{cor:H(q)} that $H((q_{s,t})_{t\geq0})$ 
  is dense in $C$, the connected 
  component  of $\real\setminus\event{\sigma=0}$ 
  containing $\alpha$.  In Corollary \ref{cor:D' totvar} we also show that almost surely 
  there is a non-empty open interval $I\subset C$ such that the total variation 
  $(D_s^{-1})'$ is infinite on any non-empty subinterval $J\subset I$. 

  Recall that by Corollary \ref{cor:rand F}, if $q_{s,t}$ and $D_s^{-1}(q_{s,t})$ are 
  semimartingales, then each point in $H(q_{s,\cdot})$ must have a neighborhood which may
  depend on $\F_s$ and on which 
  the total variation of $(D_s^{-1})'$ is finite. 
  Points in $I$ do not have this property, and
  as $H(q_{s,\cdot})$ is dense in $C$, 
  $D_s^{-1}$ 
  can not satisfy the necessary condition of Corollary \ref{cor:rand F}. We can conclude
  that $q_{s+t}=D_s^{-1}(q_{s,t})$ can not be a semimartingale in 
  $(\F_{s+t})_{t\geq0}$, which would follow from our starting assumption.
\end{proof}
 
\subsection{Continuity of \texorpdfstring{$(D'_t(x))_{t\geq0,x\in\real}$}{D'}}

First we want to show the continuity of $(t,x)\mapsto D'_t(x)$ on 
$[0,\infty)\times (\real\setminus\event{\sigma=0})$. Recall from \cite{drift:2009} that
$x\mapsto D_t(x)$ is absolutely continuous and its derivative can be written as
\begin{displaymath}
  D'_t(x)=\exp{Z_t(x)-\tfrac12\br{Z(x)}_t},
  \quad\text{where}\quad 
  Z_t(x)=\int_0^t\sigma'(D_s(x))dB_s.
\end{displaymath}
When $\abs{\sigma'}\leq L$ this yields the estimate
\begin{equation}
  \E{\abs{D'_t(x)}^p}\leq \exp{\tfrac12L^2 p(p-1)t}.
\end{equation}

\begin{lemma}\label{lem:cont of D'}
  Suppose that $\sigma$ is a convex (or concave) Lipschitz function, $B$ is a 
  Brownian motion and $(D,B)$ satisfies \eqref{eq:D (x)}.
  Then
  \begin{equation}\label{eq:Z U def}
    Z_t (x)=\int_0^t\sigma' (D_s (x))dB_s,
    \quad\text{and}\quad
    [Z(x)]_t=\int_0^t(\sigma' (D_s (x)))^2ds
  \end{equation}
  have modifications which are continuous functions of $(t,x)$ on
  $[0,\infty)\times(\real\setminus\event{\sigma=0})$  almost surely. 
\end{lemma}
\begin{remark} Even in one of  the simplest cases 
  the continuity of $D'_t(x)$ on $[0,\infty)\times\real$ may fail.
  Consider $\sigma(x)=\abs{x}$, then $Z_t(x)=\sign(x)B_t$,
  which is continuous on $[0,\infty)\times(\real\setminus\smallset{0})$
  but discontinuous at $(t,0)$ whenever $B_t\neq 0$, that is almost surely for $t>0$.
\end{remark}
\begin{proof}

  We need to show the continuity of $Z$ on $[0,\infty)\times C$, where
  $C$ is a connected component of $\real\setminus\event{\sigma=0}$.
  This is obtained via the usual combination of the
  Burkholder-Davis-Gundy inequality and Kolmogorov's lemma, see for
  example \cite[Chapter VI, proof of Theorem (1.7)]{revuz-yor}. That
  is we use that
  \begin{displaymath}
    \E{\sup_{s\leq t}\abs*{Z_s (y)-Z_s (x)}^{2p}}\leq c_p\E{\br{Z (y)-Z(x)}_t^p}
    \quad\text{for $x,y\in C$},
  \end{displaymath}
  and show that this can be further bounded by $c\abs*{y-x}^p$.
  We can assume without loss of generality that $\sigma'$ is increasing 
  (for the other case replace both $\sigma,B$ with $-\sigma,-B$).
  By assumption $\sigma'$ is bounded, say $\abs{\sigma'}\leq L$, 
  so  for $x<y$ (hence $D_s(x) \leq D_s(y)$)  we have that
  \begin{align}
      \nonumber
      \br{Z (y)-Z(x)}_t
      &=
      \int_0^t \zjel*{\sigma'(D_s(y))-\sigma'(D_s(x))}^2ds\\
      &\leq 2L \int_0^t \sigma'(D_s(y))-\sigma'(D_s(x))ds.
      \label{eq:q delta Z}    
  \end{align}

  The right hand side of \eqref{eq:q delta Z} can be related to the
  process $Y_t (x)= h (D_t (x))$, where $h:C\to\real$ is such that
  $h'=1/\sigma$. Sometimes $h (D_t (x))$ is referred to as the
  Lamperti representation of the process. 
  Simple calculus with the \ito/ formula yields that 
  \begin{displaymath}
    dY_t (x)=dB_t-\frac12\sigma' (D_t (x))dt,\quad Y_0 (x)=h (x).
  \end{displaymath}
  From this we get that
  \begin{displaymath}
    \br{Z (y)-Z(x)}_t\leq 4L \zjel*{h(y)-Y_t (y)- (h(x)-Y_t (x))}. 
  \end{displaymath}
  Then
  \begin{displaymath}
    Y_t (y)-Y_t (x)=\int_x^y h' (D_t (z)) D'_t(z)
    dz
  \end{displaymath}
  and
  \begin{displaymath}
    \zjel*{Y_t (y)-Y_t (x)}^p\leq (y-x)^p  \frac{1}{y-x} \int_{x}^y
    \abs*{h' (D_t (z))}^p 
    \abs{D'_t(z)}^{p}
    dz.
  \end{displaymath}
  To finish the proof let $a<b$  such that $[a,b]\subset C$, 
  $0<\delta<\min_{x\in[a,b]}\sigma(x)$ and
  \begin{displaymath}
    \tau_{\delta}=\tau=\inf\set*{t>0}{\text{$\sigma(D_t (z))<\delta$ for some $z\in[a,b]$}}.
  \end{displaymath}
  With these choices 
  we get that
  \begin{displaymath}
    \E{\zjel*{Y_{t\wedge\tau} (y)-Y_{t\wedge\tau} (x)}^p}
    \leq (y-x)^p \frac1{\delta^p}
    e^{\frac12 L^2 (p(p-1))t} \quad \text{if $a<x<y<b$},
  \end{displaymath}
  and 
  \begin{displaymath}
    \abs{h(y)-h(x)}^p\leq (y-x)^p\frac1{\delta^p}.
  \end{displaymath}
  
  So, we can conclude that
  \begin{displaymath}
    \E{\br{Z (y)-Z(x)}^p_{t\wedge\tau}}\leq c_p \abs*{y-x}^p,
  \end{displaymath}
  which shows that $(t,x)\mapsto Z_t (x)$ is continuous 
  on $[0,\tau)\times (a,b)$ almost surely. Letting $\delta\to0$ we have 
  that $\tau_\delta\to\infty$, so $(t,x)\mapsto Z_t (x)$ has a continuous 
  modification on $[0,\infty)\times(a,b)$ whenever $[a,b]\subset C$.
  As $C$ is open interval, this also means that there is a continuous 
  modification on $[0,\infty)\times C$.

  For $\br{Z(x)}_t$ note that by the Schwarz inequality
  \begin{displaymath}
    \zjel*{\br{Z(y)}^{1/2}_t-\br{Z(x)}^{1/2}_t}^2\leq \br{Z(x)-Z(y)}_t,
  \end{displaymath}
  hence the argument above also proves the continuity of $\br{Z(x)}^{1/2}_t$ 
  and therefore the continuity of $(x,t)\mapsto \br{Z(x)}_t$ on 
  $[0,\infty)\times\real\setminus\event{\sigma=0}$.
\end{proof}

\subsection{Decomposition of \texorpdfstring{$(q_t)_{t\geq0}$}{q\_t}}
\begin{lemma}\label{lem:Q of q}
  Suppose that $\sigma$ is Lipschitz continuous, $B$ is a Brownian motion,
  $(D,B)$ satisfies \eqref{eq:D (x)},
  and $(x,t)\mapsto D'_t(x)$ is continuous on 
  $[0,\infty)\times\real\setminus\event{\sigma=0}$ almost surely.
  For $\alpha\in\real$ let $q_t=D^{-1}_t(\alpha)$. Then
  \begin{displaymath}
    A_t = q_t+\int_0^t (D_s^{-1})'(\alpha)\sigma(\alpha) dB_s  
  \end{displaymath}
  is a process of zero energy, that is,  the quadratic variation of $A$ exists and  $\br{A}\equiv 0$.
\end{lemma}
\begin{proof}
  When $\sigma(\alpha)=0$, then $D_t(\alpha)=\alpha$  and $q_t=A_t=\alpha$ for all $t\geq0$, 
  so the claim is obvious for this case. Hence we may assume that $\sigma(\alpha)\neq 0$.

  For $s<t$
  \begin{align*}
    D_{t}(q_{t})-D_{t}(q_s)&=D_s(q_s)-D_{t}(q_s)
      =-\int_s^{t}\sigma(D_{u}(q_s))dB_u.
  \end{align*}
  From this
  \begin{align*}
    A_t-A_s
    &=q_t-q_s+\int_s^t (D^{-1}_u)'(\alpha)\sigma(\alpha)dB_u\\
    &=
    \int_s^t (D^{-1}_u)'(\alpha)\sigma(\alpha)dB_u-R(s,t)\int_s^t \sigma(D_u(q_s))dB_u,
  \end{align*}
  where by the mean value theorem $R(s,t)$ has the following form
  \begin{displaymath}
    R(s,t)=
    \begin{cases}
      \frac{q_t-q_s}{D_t(q_t)-D_t(q_s)}&q_t\neq q_s\\
      \frac1{D_t'(q_t)}&q_t=q_s      
    \end{cases}
    =\frac{1}{D_t'(\theta q_t+(1-\theta)q_s)},  
    \quad\text{for some $\theta\in(0,1)$}.
  \end{displaymath}
  The only idea in the calculation is that when $s$ is close to $t$ then 
  $\sigma(D_u(q_s))$ is close to $\sigma(\alpha)$ and $R(s,t)$ is approximately
  $(D_s^{-1})'(\alpha)$.

  Let $\pi=\smallset{t_0=0<t_1<\cdots<t_r=t}$ be a subdivision of $[0,t]$. Then
  \begin{displaymath}
    \sum_i (A_{t_{i+1}}-A_{t_i})^2 \leq 
    3\zjel*{Q_1(\pi)+Q_2(\pi)+Q_3(\pi)}
  \end{displaymath}
  where
  \begin{align*}
    Q_1(\pi) &= \sum_{i} \zjel*{\int_{t_i}^{t_{i+1}}
                \zjel*{(D^{-1}_u)'(\alpha)-(D_{t_i}^{-1})'(\alpha)}  \sigma(\alpha)dB_u}^2, \\
    Q_2(\pi) &= \sum_{i} \zjel*{R(t_i,t_{i+1})
                \int_{t_i}^{t_{i+1}} \sigma(\alpha)-\sigma(D_u(q_{t_i}))dB_u}^2, \\
    Q_3(\pi) &= \sum_{i} \zjel*{((D_{t_i}^{-1})'(\alpha)-R(t_i,t_{i+1}))
                \int_{t_i}^{t_{i+1}} \sigma(\alpha)dB_u}^2.
  \end{align*}
  For $Q_3$ observe that since $(t,x)\mapsto (D_t^{-1})'(x)$ is continuous, we have 
  that
  \begin{displaymath}
    \max_{i} \abs*{(R(t_i,t_{i+1})-(D_{t_i}^{-1})'(\alpha)}\to 0,
    \text{almost surely as $\max_{i}(t_{i+1}-t_i)\to0$},
  \end{displaymath}
  so $Q_3(\pi)\pto0$ as the mesh of the partition $\pi$ goes to zero.

  For $Q_1$, $Q_2$ we use the next observation
  \begin{proposition}\label{prop:Q to 0} 
    For a subdivision $\pi=\smallset{t_0=0<t_1<\cdots<t_r=t}$ of $[0,t]$ let 
    $\pi(u)=t_i$ if $u\in[0,t)$ and $t_i\leq u<t_{i+1}$, that is,
    $\pi(u)$ is the starting point of the subinterval containing $u$. 
    Denote $\mesh(\pi)=\max_i (t_{i+1}-t_i)$.
    Suppose that $(\phi(u,v))_{u,v\geq0}$ is a two parameter process 
    such that $\phi(u,v)\I{u\geq v}$ is integrable with respect to the Brownian motion $B$
    for each fixed $v$.    
    If 
    \begin{equation}\label{eq:Q to 0 cond}
      \int_0^t \phi^2(u,\pi(u))du\pto0, \quad\text{as $\mesh(\pi)\to0$},
    \end{equation}
    then
    \begin{displaymath}
      \sum_i \zjel*{\int_{t_i}^{t_{i+1}} \phi(u,t_i) dB_u}^2\pto 0,\quad\text{as $\mesh(\pi)\to0$}.
    \end{displaymath}
  \end{proposition}
  For $Q_1(\pi)$ we can use 
  \begin{displaymath}
    \phi(u,v)=(D^{-1}_u)'(\alpha)-(D_{v}^{-1})'(\alpha)   
  \end{displaymath}
  and note that by the continuity of $(x,u)\mapsto D^{-1}_u(x)$ 
  condition \eqref{eq:Q to 0 cond} holds and $Q_1(\pi)\pto0$ as $\mesh(\pi)\to0$.

  For $Q_2$ we consider
  \begin{displaymath}
    \phi(u,v)= \sigma(D_u(q_{v}))-\sigma(\alpha)= \sigma(D_u(q_{v}))-\sigma(D_v(q_{v}))
  \end{displaymath}
  and now we can refer to the continuity of $(x,u)\mapsto D_u(x)$ and $q$ to get that 
  condition \eqref{eq:Q to 0 cond} holds. To see that $Q_2(\pi)\pto0$ we only need to add 
  that 
  \begin{displaymath}
    \max_{i} R(t_i,t_{i+1})\to \max_{s\leq t } (D^{-1}_s)'(\alpha)<\infty,
    \quad\text{almost surely as $\mesh(\pi)\to0$}. \qedhere
  \end{displaymath}
\end{proof}

\begin{proof}[Proof of Proposition \ref{prop:Q to 0}]
  For  a subdivision $\pi=\smallset{t_0=0<t_1<\cdots<t_r=t}$ of $[0,t]$,
  denote by $(Y(s,\pi))_{s\geq0}$ the following process
  \begin{displaymath}
    Y(s)=Y(s,\pi) = \int_0^s \phi(u,\pi(u))dB_u.
  \end{displaymath}
  Note that from $\int_0^t \phi^2(u,\pi(u))du\to0$ it  follows that
  \begin{equation}\label{eq:sup Y}
    \sup_{s\leq t}\abs{Y(s,\pi)}\pto0\quad\text{as $\mesh(\pi)\to0$}.
  \end{equation}

  By \ito/'s formula
  \begin{displaymath}
    \zjel*{\int_{t_i}^{t_{i+1}} \phi(u,t_i)dB_u}^2 = 
    2\int_{t_i}^{t_{i+1}} (Y(u)-Y(t_i))\phi(u,t_i)dB_u 
    + \int_{t_i}^{t_{i+1}} \phi^2 (u,t_i) du.
  \end{displaymath}
  Hence
  \begin{displaymath}
    \sum_i \zjel*{\int_{t_i}^{t_{i+1}} \phi(u,t_i)dB_u}^2 =
    2\int_0^t (Y(u)-Y(\pi(u)))\phi(u,\pi(u))dB_u 
    + \int_0^t \phi^2 (u,\pi(u)) du.
  \end{displaymath}
  By assumption the time integral tends to zero in probability as $\mesh(\pi)\to0$.
  For the stochastic integral it is enough to show that
  \begin{displaymath}
    \int_0^t (Y(u)-Y(\pi(u)))^2\phi^2(u,\pi(u))du\pto0.
  \end{displaymath}
  Here
  \begin{displaymath}
    \int_0^t (Y(u)-Y(\pi(u)))^2\phi^2(u,\pi(u))du\leq 
    4\sup_{u\leq t} Y^2(u) \int_0^t \phi^2(u,\pi(u))du.
  \end{displaymath}
  Although it is suppressed in the notation both factors in the upper bound depend on $\pi$. 
  Nevertheless, it is easy to see that both factor tends to zero in probability and so does 
  their product. For the second factor it is the main assumption, while for the first one it 
  was already observed in  \eqref{eq:sup Y}.
\end{proof}

The next claim follows from Lemma \ref{lem:Q of q}.
\begin{corollary}\label{cor:H(q)}
  Suppose that the assumptions of Lemma \ref{lem:Q of q} hold and $\sigma(\alpha)\neq0$.
  Let $C$ be the connected component of 
  $\real\setminus\event{\sigma=0}$ containing $\alpha$. Then
  \begin{displaymath}
    \P{\int\nolimits_0^\infty \I{q_s\in(a,b)}d\br{q}_s>0}>0
    \quad\text{for any $(a,b)\subset C$}.
  \end{displaymath}
  Especially, if $q$ is  a semimartingale, then $H(q)$  is dense in $C$.
\end{corollary}
\begin{proof}
  By Lemma \ref{lem:Q of q} $q$ has quadratic variation and
  \begin{displaymath}
    \int_0^\infty \I{q_s\in(a,b)}d\br{q}_s=
    \int_0^\infty \I{D_s(a)<\alpha<D_s(b)} ((D_s^{-1})'(\alpha)\sigma(\alpha))^2 ds.
  \end{displaymath}
  Since $((D_s^{-1})'(\alpha)\sigma(\alpha))^2>0$, it is enough to show that the 
  continuous process $(D_t(x))_{t\geq0}$ hits $\alpha$ with positive probability 
  for any $x$ in $C$. Indeed, the set $\set{s>0}{D_s(a)<\alpha<D_s(b)}$ is open, so 
  if not empty, then it has positive Lebesgue measure.
  
  This is an easy martingale argument. Let $c<\alpha,x<d$ be such that $[c,d]\subset C$. Then
  $\inf\sigma|_{[c,d]}>0$ implies that $\tau=\inf\set*{t\geq0}{D_t(x)\notin[c,d]}$ is 
  almost surely finite, and since $\E{D_\tau(x)}=x$, both $\P{D_\tau(x)=c}$ and $\P{D_\tau(x)=d}$ 
  are positive, which implies by the continuity of the sample paths of $D(x)$ that 
  $D(x)$ hits $\alpha$ with positive probability.

  Finally, suppose that $q$ is a semimartingale. Then its local time exists, and 
  for any fixed $t$ the mapping $x\mapsto L^x_t(q)$ is \cadlag/.
  From the first part of the statement, for $(a,b)\subset C$  and $t$ large enough
  $\int_{(a,b)} L^x_t(q)dx>0$ with positive probability. This yields that for some $\eps>0$
  and $(a',b')\subset (a,b)$  we have $\inf_{x\in (a',b')} L^x_t(q)>\eps$ with positive 
  probability and $(a',b')\subset H(q)$.
\end{proof}

\subsection{Quadratic variation of \texorpdfstring{$D_t'$}{D'\_t} in the space variable}
\begin{lemma}\label{lem:q Z}
  Suppose that $\sigma$ is Lipschitz continuous and is the difference of convex functions,
  $B$ is a Brownian motion,
  $(D,B)$ satisfies \eqref{eq:D (x)},
  and 
  \begin{displaymath}
    (x,t)\mapsto Z_t(x)=\int_0^t \sigma'(D_s(x))dB_s  
  \end{displaymath}
  is continuous on 
  $[0,\infty)\times\real\setminus\event{\sigma=0}$ almost surely.
Then
the quadratic variation processes of the random functions $x\mapsto
Z_t (x)$ and $x\mapsto U_t(x)=\br{Z(x)}_t$ exist,
\begin{displaymath}
  \br{Z_t} (b)-\br{Z_t} (a)= \int_0^t \sum_{z\in [D_s(a),D_s(b))} (\Delta\sigma'(z))^2 ds 
  \quad\text{and}\quad
  \br{U_t}\equiv 0,
\end{displaymath}
where $\sigma'$ denotes the left hand side derivative of $\sigma$.
\end{lemma}
\begin{proof}
\begin{align*}
  (Z_t (y)-Z_t (x))^2 
  &=2\int_0^t(Z_s (y)-Z_s (x)) (\sigma' (D_s (y))-\sigma' (D_s (x)))dB_s\\
  &\hphantom{{}={}}
  +\int_0^t(\sigma' (D_s (y))-\sigma' (D_s (x)))^2ds.
\end{align*}
For a subdivision $a=x_0<x_1<\dots<x_r=b$
\begin{multline}
  \sum_i(Z_t (x_{i+1})-Z_t (x_i))^2
  \\
  \begin{aligned}
    &=2\int_0^t\sum_i(Z_s (x_{i+1})-Z_s (x_i))
    (\sigma' (D_s (x_{i+1}))-\sigma' (D_s (x_i)))dB_s
  \\
  &\hphantom{{}={}}
  +\int_0^t\sum_i(\sigma' (D_s (x_{i+1}))-\sigma' (D_s (x_i)))^2ds. 
  \end{aligned}
  \label{eq:q Z}
\end{multline}
As $\sigma$ is the difference of convex functions, its left derivative 
exists everywhere and is of finite total variation on compact intervals. 
Denote by $\mu$ the total variation measure of $\sigma'$. 
Then the integrand in the stochastic integral is bounded by
\begin{displaymath}
  \sup_{s\leq t,\,i} \abs*{Z_s (x_{i+1})-Z_s (x_i)} \mu
  ([\min_{s\leq t}D_s (a),\max_{s\leq t}D_s (b)]).
\end{displaymath}
This upper bound goes to zero almost surely by the continuity of $(s,x)\mapsto
Z_s (x)$, hence the integral with respect to the Brownian motion
goes to zero in probability as the mesh of the subdivision goes to zero. 

For the second integral, note that
\begin{align*}
  \sum_i(\sigma' (D_s (x_{i+1}))-\sigma' (D_s (x_i)))^2
  &\to \sum_{z\in [D_s(a),D_s(b))} (\Delta \sigma'(z))^2
  \intertext{and}                                                      
  \sum_i(\sigma' (D_s (x_{i+1}))-\sigma' (D_s (x_i)))^2
  &\leq (\mu ([\min_{s\leq t}D_s (a),\max_{s\leq t}D_s (b)))^2.
\end{align*}
So by the dominated convergence theorem the time integral in 
\eqref{eq:q Z}  almost surely goes to
\begin{displaymath}
  \int_0^t \sum_{z\in[D_s (a),D_s (b))}(\Delta \sigma'(z))^2 ds
\end{displaymath}
as the mesh of the subdivision goes to zero and the
formula for the quadratic variation of $x\mapsto Z_t (x)$ is proved.

The computation for the $U_t$ is similar
\begin{align*}
  (U_t (y)-U_t (x))^2
  &=2\int_0^t (U_s (y)-U_s(x))
    \zjel*{(\sigma')^2 (D_s(y))-(\sigma')^2 (D_s (x))}ds\\
  &\leq
    4L \int_0^t\abs*{U_s (y)-U_s
    (x)}\mu ([D_s (x),D_s (y)))ds,
\end{align*}
where $L$ is the Lipschitz constant for $\sigma$, that is
$\sup_{r}\abs*{\sigma' (r)}$,
and
\begin{displaymath}
  \sum_i(U_t (x_{i+1})-U_t (x_i))^2\leq
  4Lt\mu ([\min_{s\leq t}D_s (a),\max_{s\leq t}D_s (b)))
  \sup_{s\leq t,\,i}\abs*{U_t (x_{i+1})-U_t (x_i)}.
\end{displaymath}
By the continuity of $(s,x)\mapsto U_s (x)$
this upper bound goes to zero almost surely  
when the mesh of the subdivision goes to zero.
\end{proof}

\begin{lemma}\label{lem:q of h(X)}
  Let $(X_z)_{z\in[a,b]}$ be a process with continuous sample paths
  whose quadratic variation $(\br{X}_z)_{z\in[a,b]}$ exists.
  If $h:[a,b]\to\real$ is continuously differentiable, then 
  the  quadratic variation of $h(X)$ exists and 
  \begin{equation}
    \label{eq:Q h(X)}
    \br{h(X)}(b)-\br{h(X)}(a)=\int_a^b \zjel*{h' (X_z)}^2\br{X}(dz).
  \end{equation}
\end{lemma}
\begin{proof} 
  When $h$ is linear or continuous and piecewise linear, then \eqref{eq:Q h(X)} 
  follows from the definition of the quadratic variation.
  
  For a subdivision $\pi=\smallset{t_0=a<t_1<\dots<t_r=b}$ of $[a,b]$ 
  and for a process $Y$ denote 
  \begin{displaymath}
    Q^2_\pi(Y)=\sum_{i} (Y_{t_{i+1}}-Y_{t_i})^2.
  \end{displaymath}
  Note that $Q_\pi(Y)$ is the Euclidean norm of the vector 
  $(Y_{t_1}-Y_{t_0},\dots,Y_{t_r}-Y_{t_{r-1}})$, so by the triangle inequality,
  for any processes $Y,Z$ and subdivision $\pi$, we have
  \begin{displaymath}
    \abs{Q_\pi(Y)-Q_\pi(Z)}\leq Q_\pi(Y-Z).  
  \end{displaymath}
  Similarly, if $\tQ^2(h)=\int_a^b (h'(Z_z))^2d\br{Z}(dz)$ denotes the right 
  hand side of \eqref{eq:Q h(X)}, then
  \begin{displaymath}
    \abs{\tQ(h_1)-\tQ(h_2)}\leq \tQ(h_1-h_2).
  \end{displaymath}

  If $h$ is Lipschitz continuous  with Lipschitz constant $L$, then
  \begin{displaymath}
    (h(X_t)-h(X_s))^2\leq L^2(X_t-X_s)^2,
  \end{displaymath}
  so for any subdivision $\pi$ 
  \begin{displaymath}
    Q_{\pi}(h(X))\leq L Q_{\pi}(X).
  \end{displaymath}
  When $h$ is $C^1$ and $\eta>0$, then there is a continuous piecewise linear function $h_\eta$ 
  such that $h-h_\eta$ is Lipschitz continuous with Lipschitz constant $\eta$. 
  Then 
  \begin{align*}
    &\abs*{Q_\pi(h(X))- \tQ(h)} \\ 
    &\hphantom{Q^{1/2}}\leq
    \abs*{Q_\pi(h(X))-Q_\pi(h_\eta(X))} 
    +\abs*{Q_\pi(h_\eta(X))-\tQ(h_\eta)} 
    +\abs{\tQ(h_\eta)-\tQ(h)}\\
    &\hphantom{Q^{1/2}}\leq Q_\pi((h-h_\eta)(X))
    +\abs*{Q_\pi(h_\eta(X))-\tQ(h_\eta)}
    +\tQ((h-h_\eta))
    \\
    &\hphantom{Q^{1/2}}\leq \eta Q_\pi(X)+
    \abs*{Q_\pi(h_\eta(X))-\tQ(h_\eta)}+\eta(\br{X}_b-\br{X}_a)^{1/2}.
  \end{align*}
  Now, let $(\pi_n)$ be a sequence of subdivisions of $[a,b]$ with $\mesh(\pi_n)\to0$. 
  From the previous estimation, using that the middle term and 
  $Q_{\pi_n}(X)-(\br{X}_b-\br{X}_a)^{1/2}$ goes to zero in probability,
  we get that
  \begin{displaymath}
    \limsup_{n}\P{\abs*{Q_{\pi_n}(h(X))- \tQ(h)}>2\eps}
    \leq \inf_{\eta>0} \P{2\eta(\br{X}_b-\br{X}_a)^{1/2}>\eps}=0.
  \end{displaymath}
  So $Q_{\pi_n}(h(X))\pto \tQ(h)$ and $Q^2_{\pi_n}(h(X))\pto \tQ^2(h)$, 
  which proves the claim. 
\end{proof}

\begin{corollary}
  Under the assumptions of Lemma \ref{lem:q Z} the quadratic variation 
  of $x\mapsto D'_t(x)$ exists and is given by the formula
  \begin{displaymath}
    \br{D'_t}(b)-\br{D'_t}(a) = \int_{a}^b 4(D'_t(x))^2 \br{Z_t}(dx),
  \end{displaymath}
  where
  \begin{displaymath}    
    \br{Z_t}(x)-\br{Z_t}(a)=
    \int_0^t \sum_{z\in [D_s(a),D_s(x))}(\Delta\sigma' (z))^2 ds.
  \end{displaymath}
\end{corollary}

\begin{corollary}\label{cor:D' totvar prep}
  Suppose that the assumptions of Lemma \ref{lem:q Z} hold.
  Denote by $S$ the set of discontinuity points of $\sigma'$ and let $C$ be 
  a connected component of $\real\setminus\event{\sigma=0}$.

  On an almost sure event the following holds: 
  if  
  \begin{displaymath}
    (a,b)\subset C\quad\text{and}\quad
    (\min_{s\leq t} D_s(a),\max_{s\leq t}D_s(b))\cap S\neq\emptyset,
  \end{displaymath}
  then the total variation of $D'_t$ over the interval $[a,b]$ is infinite.
\end{corollary}

\begin{corollary}\label{cor:D' totvar}
  Suppose that the assumptions of Lemma \ref{lem:q Z} hold.
  Denote by $S$ the set of discontinuity points of $\sigma'$ and let $C$ be 
  a connected component of $\real\setminus\event{\sigma=0}$.
  Suppose that $S\cap C\neq\emptyset$

  Then almost surely there exists a non-empty open interval $I\subset C$ 
  such that the total variation of $x\mapsto D'_s(x)$ is infinite on any non-empty 
  subinterval of $I$.

  The same is true for $(D_s^{-1})'$.
\end{corollary}
\begin{proof}
  Let $z\in S\cap C$.  Since $\sigma(z)\neq0$, we have that 
  $\min_{t\leq s} D_t(z)<z<\max_{t\leq s} D_t(z)$. By the continuity of $(t,x)\mapsto D_t(x)$
  there is $a<z<b$, such that $\min_{t\leq s}D_s(b)<z<\max_{t\leq s} D_t(a)$. 
  As $D_t(x)$ is increasing in $x$ from this we have that for any $a<c<d<b$
  \begin{displaymath}
    \min_{t\leq s} D_t(c)<\min_{t\leq s} D_t(b)<z<\max_{t\leq s} D_t(a)<\max_{t\leq s} D_t(d).
  \end{displaymath}
  So, by Corollary \ref{cor:D' totvar prep} the total variation of $D'_s$ on $(c,d)$ is infinite.

  Finally, since $(D_s^{-1})' = 1/ (D_s'\circ D_s^{-1})$ and $D_s$ maps $C$ onto $C$, the 
  image of $I$ under $D_s$ is a subinterval of $C$ such that the total variation of 
  $(D^{-1}_s)'$ is infinite on any of its non-empty subintervals.
\end{proof}

\section{Some technical results}

\begin{lemma}\label{lem:ell}
  Let $X$ be a one dimensional continuous semimartingale
  and for $x\in\real$, $\delta>0$ put
  \begin{align*}
    \sigma^{x,\delta}_0 & = 0,\\
    \sigma^{x,\delta}_{2k+1}
    &=\inf\set*{t\geq\sigma^{x,\delta}_{2k}}{X_t=x},\quad k\geq0,\\ 
    \sigma^{x,\delta}_{2k+2}
    &=\inf\set*{t\geq\sigma^{x,\delta}_{2k+1}}{\abs*{X_t-x}=\delta},\quad k\geq0,\\
    \ell^\delta(x,t)&=\inf\set*{k}{\sigma^{x,\delta}_{2k+1}>t}=\sum_{k\geq0}\I{\sigma_{2k+1}\leq t}.
  \end{align*}

  Suppose that $(\rho_K)$ is a sequence of stopping times tending to
  infinity almost surely as $K\to\infty$ and $x\in H(X)$. Then there
  is a $K>0$, a non-empty open interval $I$ containing $x$ and $\delta_0>0$ 
  such that
  \begin{displaymath}
    \inf_{\delta\in (0,\delta_0)}\inf_{y \in I}\delta\E{\ell^{\delta} (y,\rho_K)}>0.
  \end{displaymath}
\end{lemma}
\begin{remark}
Note that 
the definition of $\ell^\delta$ in this claim is slightly
different from that of used in Theorem \ref{thm:BM} and \ref{thm:X}.
In this lemma it is defined for $x\in\real$, while in the 
formulation with the embedded random walk it was only defined for $x\in\delta\Z$.
\end{remark}

We will denote by $X=X_0+M+A$ the semimartingale decomposition of $X$, 
that is $M$ is local martingale starting from zero and $A$ is process of finite variation.
Also, $V$ is used for the total variation process of $A$.

The proof of Lemma \ref{lem:ell} is based on the following estimation.
\begin{proposition}\label{prop:L-ell}
  With the notation introduced above let 
  \begin{displaymath}
    \ell_{+}^\delta(x,t)=\sum_{k\geq 0}\I{X_{\sigma_{2k}>x}}\I{\sigma_{2k+1}\leq t}.
  \end{displaymath}
  If $\rho$ is a stopping time such that $\br{M}_\rho$ and $V_\rho$ are integrable, 
  then
  \begin{displaymath}
    \limsup_{\delta\to0}\sup_{x\in[a,b)}\E{\frac12 L^x_\rho-\delta \ell_{+}^\delta(x,\rho)}
    \leq
    \sup_{x\in(a,b]}\E{\int_0^\rho \I{X_s=x}dV_s}.
  \end{displaymath}
  where $(L^x_t)_{x\in\real,t\geq0}$ is the family of local times for $X$.
\end{proposition}
\begin{proof}[Proof of Lemma \ref{lem:ell}]
  $\rho_K$ is a sequence of stopping time tending to infinity almost surely. 
  Since  $\ell^\delta(x,t)$ is non-decreasing in $t$, we can clearly assume
  $\rho_K$ is bounded, and $\br{M}_{\rho_K}$, $V_{\rho_K}$ are integrable for each $K$.
  For $x\in H(X)$, $\E{L^x_{\rho_K}}$, $\E{L^{x-}_{\rho_K}}>0$ for some $K$ sufficiently large.

  By our assumptions on $\rho_K$ the function $\E{L^x_\rho}$ is right continuous, 
  so there is $b>x$ such that $\inf_{y\in[z,b)}\E{L^y_\rho}=c\geq \frac12 \E{L^x_\rho}>0$.
  The formula $\mu(C)=\E{\int_0^\rho \I{X_s\in C}dV_s}$ defines a finite measure 
  on the Borel $\sigma$-algebra of the real line, whence there are only finitely many $y$-s,
  such that $\mu({\smallset{y}})>c/4$, and taking a smaller $b$ if necessary it is possible 
  to achieve that $\sup_{y\in(x,b]}\mu(\smallset{y})<c/4$. With this choice
  \begin{displaymath}
    \liminf_{\delta\to 0}\inf_{y\in[x,b)}\E{\delta\ell^{\delta}(y,\rho)}
    \geq\liminf_{\delta\to 0}\inf_{y\in[x,b)}\E{\delta\ell^{\delta}_{+}(y,\rho)}
    > \tfrac14c\geq \tfrac18 \E{L^x_\rho}>0
  \end{displaymath}
  by  Proposition \ref{prop:L-ell}.

  Similar analysis applies to the left hand side, as $L^{-x}_t(-X)=L^{x-}_t(X)$. So from 
  the previous argument we get $a<x$, such that 
  \begin{displaymath}
    \liminf_{\delta\to 0}\inf_{y\in(a,x]}\E{\delta\ell^{\delta}(y,\rho)} > \tfrac18\E{L^{x-}_\rho}>0
  \end{displaymath}

  Then with this $K$ we have a $\delta_0$,  such that for $a<x<b$ 
  \begin{displaymath}
    \inf_{y\in(a,b),\,0<\delta<\delta_0} \delta\E{\ell^\delta(y,\rho_K)} >0,
  \end{displaymath}
  which completes the proof.
\end{proof}

\begin{proof}[Proof of Proposition \ref{prop:L-ell}]
  Let 
  \begin{displaymath}
    \phi^{\delta}_s
    =\sum_{k \geq 0}
    \I{\sigma_{2k}\leq s \leq\sigma_{2k+1}}
    \I{X_{\sigma_{2k}}>x}.
  \end{displaymath}
  Then
  \begin{align*}
    \int_0^\rho \phi^{\delta}_s dX_s
    &=
    \sum_{k\geq0}\I{X_{\sigma_{2k}}>x}\zjel{X_{\sigma_{2k+1}\wedge\rho}-X_{\sigma_{2k}\wedge\rho}} 
  \end{align*}
  and 
  from the definition of the stopping times $(\sigma_k)_{k\geq0}$,
  we obtain that
  \begin{displaymath}
    \abs*{\abs*{X_{\rho}-x}_+
    -\abs*{X_0-x}_{+}
    -\int_0^{\rho}\phi^{\delta}_s dX_s
    -\delta\ell^{\delta}_+(x,\rho)}\leq 2\delta.
  \end{displaymath}
  Using the Tanaka formula, we have
  \begin{displaymath}
    \frac12 L^x_{\rho}
    =\abs*{X_{\rho}-x}_+-\abs*{X_0-x}_{+}-\int_0^{\rho}\I{X_s>x}dX_s.
  \end{displaymath}
  The mean difference between $\frac12 L^x_{\rho}$ and
  $\delta\ell^{\delta}_+(x,\rho)$ (using that $s \mapsto \phi^{\delta}_s-\I{X_s>x}$ is
  bounded and $E{\br{M}_{\rho_K}}<\infty$, so the expected value of the integral w.r.t. $M$ vanishes) 
  is 
  \begin{align*}
    \E{\frac12 L^x_{\rho}-\delta\ell^{\delta}_{+} (x,\rho)}
    &\leq 2\delta + \E{\int_{0}^{\rho}\phi^{\delta}_s-\I{X_s>x}dX_s} \\
    &= 2\delta + \E{\int_{0}^{\rho_K}\phi^{\delta}_s-\I{X_s>x}dA_s}.
  \end{align*}
  Recall that  $V$ is the total variation process of $A$. We have the following estimate
  \begin{align}
    \abs*{\E{\frac12 L^x_{\rho}-\delta\ell^{\delta}_{+} (x,\rho_K)}}
    &\leq2\delta + \E{\int_{0}^{\rho}\abs*{\phi^{\delta}_s-\I{X_s>x}}dV_s} \nonumber \\
    &\leq2\delta + \E{\int_{0}^{\rho}\I{X_s\in (x,x+\delta)}dV_s}. \label{eq:diff loc time}
  \end{align}

  Let $g$ and $h$ be the following functions
  \begin{displaymath}
    g(y, \delta) = \E{\int_{0}^{\rho}\I{X_s\in (y,y+\delta)}dV_s}, \quad
    h(y)=\E{\int_0^{\rho}\I{X_s=y}dV_s}.
  \end{displaymath} 
  For an interval $[a, b)$ let $(y_n, \delta_n)_{n \geq 1}$ be a sequence 
  for which 
  \begin{displaymath}
    \lim_{n\to\infty} g(y_n, \delta_n) = 
    \limsup_{\delta\to0}\sup_{y \in [a,b)}g(y,\delta).  
  \end{displaymath}
  Clearly, we can assume (by passing to a subsequence if necessary) that 
  $\delta_n \to 0$ and $y_n$ is convergent. 
  Let $y^* = \lim_{n\to\infty}y_n$, $y^*\in[a,b]$. 
  We have the following estimation:
  \begin{displaymath}
    g(y_n, \delta_n) \leq 
    \E{\int_0^{\rho}
    \I{X_s \in (\min(y_n,y^*), \max(y_n+\delta_n, y^*+\delta_n))}dV_s}.
  \end{displaymath}
  
  If $y_n < y^*$ only for finitely many $n$, then for $n$ large enough we have 
  $y_n \geq y^*$ and 
  \begin{displaymath}
    \E{\int_0^{\rho}
    \I{X_s \in (\min(y_n,y^*), \max(y_n+\delta_n, y^*+\delta_n))}dV_s}  
    \leq 
    \E{\int_0^{\rho}
    \I{X_s \in (y^*, y^*+\delta_n)}dV_s},
  \end{displaymath}
  and by the dominated convergence theorem the right hand side tends to $0$. 
  If $y^*=a$, then this is the case.
  
  If $y_n < y^*$ for infinitely many $n$, then
  \begin{align*}
    \limsup_{\delta\to0}\sup_{y \in [a,b)}g(y,\delta) 
    &= 
    \lim_{n\to\infty} g(y_n,\delta_n) 
    \\ 
    &\leq \limsup_{n \to \infty} 
    \E{\int_0^{\rho}
    \I{X_s \in (\min(y_n,y^*), \max(y_n+\delta_n, y^*+\delta_n))}dV_s} 
    \\
    &\leq \E{\int_0^{\rho}\I{X_s = y^*}dV_s} = h(y^*).
  \end{align*}
  In this case $y^* \in (a,b]$. So in both cases we have
  \begin{displaymath}
    \limsup_{\delta \to0} \sup_{y \in [a,b)}g(y, \delta) \leq 
    \sup_{y \in (a,b]} h(y).    
  \end{displaymath}
  By \eqref{eq:diff loc time} this proves the claim.
\end{proof}

Recall that 
\begin{displaymath}
  \mu_{\delta,F}(H)=\sum_{x\in H\cap\delta \Z} \frac1\delta \abs{(\Delta^\delta F) (x)},
\end{displaymath}
where $(\Delta^\delta F)(x)=F(x+\delta)+F(x-\delta)-2F(x)$.
\begin{proposition}\label{prop:tot var} Let $F:\real\to\real$ be a
  continuous function and $I$ be an open interval.
  Then $F$ is the difference of convex 
  functions on $I$ if and only if 
  \begin{equation}\label{eq:sup mu_delta}
    \limsup_{\delta\to0}\mu_{\delta,F}(J)<\infty,
    \quad\text{for all compact intervals $J\subset I$} .
  \end{equation}
\end{proposition}
\begin{proof}
  First suppose that  \eqref{eq:sup mu_delta} holds for $I$. 
  It is enough to show that $F$ is a difference of convex functions on compact intervals 
  $J\subset I$. Using our assumption we define below two sequences of  convex functions 
  $f_n,g_n:J\to\real$ such that  $F-(f_n-g_n)$ is 
  linear on $J\cap 2^{-n}\Z$  for each $n$ and $f_{n_k}\to f$, $g_{n_k}\to g$  pointwise
  for suitable subsequence $(n_k)$ of the indices.  Then $f,g$ are convex and $F-(f-g)$ is 
  a linear function on $J$, hence $F$ is a difference of convex functions on $J$.

  Let us now define $f_n,g_n$. For $\delta>0$ and $H\in\B(J)$ let
  \begin{align*}
    \tmu_{\delta,+} (H)
    &=\sum_{x\in H\cap \delta\Z} \frac1\delta \abs{(\Delta^\delta F)(x)}_{+},
    &\tmu_{\delta,-} (H)
    &=\sum_{x\in H\cap \delta\Z}
      \frac1\delta \abs{(\Delta^\delta F)(x)}_{-}.
  \end{align*}
  $\tmu_{\delta,+}$, $\tmu_{\delta,-}$ are finite measures on the 
  Borel $\sigma$-algebra of $J$. 
  Since $\phi(x,y)=\abs{x-y}_{+}$ is 
  a convex function in $x$ for each fixed $y$, we have that
  \begin{displaymath}
    f_\delta(x) = \int_J \phi(x,y)\tmu_{\delta,+}(dy),
    \quad
    g_\delta(x) = \int_J \phi(x,y)\tmu_{\delta,-}(dy)
  \end{displaymath}
  are convex functions (on $\real$). For $x,y\in\delta \Z$, we have that
  \begin{displaymath}
    \Delta^\delta \phi(x,y)=\phi(x+\delta,y)+\phi(x-\delta,y)-2\phi(x,y)=
    \begin{cases}
      0&x\neq y\\
      \delta& x=y
    \end{cases},
  \end{displaymath}
  which yields for $x\in J\cap\delta Z$ that
  \begin{displaymath}
    \frac1\delta (\Delta^\delta f_\delta)(x)=  \mu_{\delta,+}(\smallset{x}),
    \quad
    \frac1\delta (\Delta^\delta g_\delta)(x)=  \mu_{\delta,-}(\smallset{x})
  \end{displaymath}
  and
  \begin{displaymath}
    \frac1\delta \Delta^{\delta}(F-(f_\delta-g_\delta))(x)=0,
    \quad \text{for all $x\in J\cap \delta\Z$}.
  \end{displaymath}
  We obtained that $F-(f_\delta-g_\delta)$ is linear on $J\cap \delta Z$.

  We put $f_n=f_{2^{-n}}$ and $g_n=g_{2^{-n}}$. 
  The families of  finite measures $\set{\mu_{2^{-n},+}}{n\geq1}$ and
  $\set{\mu_{2^{-n},-}}{n\geq1}$ are tight as $J$ is compact and 
  \begin{displaymath}
    \sup_{\delta\in(0,1)}\mu_{\delta,\pm}(J)
    \leq \sup_{\delta\in(0,1)}\mu_{\delta,F}(J)<\infty.
  \end{displaymath}
  Then there is 
  a subsequence $(n_k)$ of the indices such that $\mu_{n_k,+}$ and $\mu_{n_k,-}$ 
  converge weakly to some limit. As  $y\mapsto\phi(x,y)$ is bounded continuous on $J$, 
  the sequences $f_{n_k}(x)=\int \phi(x,y)\mu_{n_k,+}(dy)$ and 
  $g_{n_k}(x)=\int \phi(x,y)\mu_{n_k,-}(dy)$ 
  are pointwise convergent on $J$, which finishes the proof 
  of the suffiency of \eqref{eq:sup mu_delta}.
  
  To prove the necessity of \eqref{eq:sup mu_delta} it is enough to consider a convex $F$.
  Let $[a,b]\subset  I$ and $\delta$ so small that $(a-\delta,b+\delta)\subset I$.
  As  $F$  is convex on $I$,  we have that  $\Delta^\delta F (x)\geq 0$ on $[a,b]$
  and
  \begin{align*}
    \mu_{\delta,F}([a,b])=\sum_{x\in J\cap \delta Z}\frac1\delta (\Delta^\delta F)(x)
    &=
    \frac{F(b'+\delta)-F(b')}{\delta}-\frac{F(a')-F(a'-\delta)}{\delta}
    \\
    &\leq
    \frac{F(b+\delta)-F(b)}{\delta}-\frac{F(a)-F(a-\delta)}{\delta},
  \end{align*}
  where $a'=\min [a,b]\cap \delta \Z$ and $b'=\max [a,b]\cap \delta Z$. 
  In the first step we used that we have a telescoping sum, while in the last step we 
  used that $a\leq a'\leq b'\leq b$ and the divided difference 
  $x\mapsto \frac1\delta (F(x+\delta)-F(x))$ is increasing in $x$. Now letting $\delta\to 0$ 
  gives that
  \begin{displaymath}
    \limsup_{\delta\to0}\mu_{\delta,F}([a,b])\leq F'_{+}(b)-F'_{-}(a). 
  \end{displaymath}
  This upper bound is finite since $F$ is convex on $I$ and $[a,b]\subset\interior I$.
\end{proof}

\bibliography{main_bib}
\bibliographystyle{elsarticle-harv}

\end{document}